\documentclass[smallextended,envcountsect]{svjour3}
\usepackage{graphicx,color,cite}
\usepackage{mathrsfs}
\usepackage{amsmath, amscd, amsfonts, amssymb, graphicx, color}
\usepackage[bookmarksnumbered, colorlinks, plainpages]{hyperref}
\usepackage{mathptmx}   
\usepackage[colorinlistoftodos]{todonotes}
\usepackage[left=4cm, right=4cm, top=3.5cm, bottom=3.5cm]{geometry}

\usepackage{latexsym}

\def\tcv{ \textcolor{violet}}
\definecolor{vert_10}{RGB}{21,106,47}
\if{
\newtheorem{theorem}{Theorem}[section]
\newtheorem{definition}{Definition}[section]
\newtheorem{lemma}{Lemma}[section]

\newtheorem{remark}{Remark}[section]
\newtheorem{example}{Example}[section]

\newtheorem{corollary}{Corollary}[section]

}\fi
\newcommand{\R}{\mathbb R}

\def\dom{\mathop{\rm dom\,}}
\def\epi{\mathop{\rm Epi\,}}

\def\cone{\mathop{\rm cone\,}}

\renewcommand\theenumi{\roman{enumi}}

\newcounter{mycount}

\smartqed

\usepackage{etex}

\usepackage{mathtools}

\usepackage{enumitem}
\renewcommand\theenumi{(\roman{enumi})}
\renewcommand\theenumii{(\alph{enumii})}

\usepackage{csquotes}

\usepackage{todonotes}

\makeatletter
\let\orgdescriptionlabel\descriptionlabel
\renewcommand*{\descriptionlabel}[1]{
 \let\orglabel\label
 \let\label\@gobble
 \phantomsection
 \edef\@currentlabel{#1}
 \let\label\orglabel
 \orgdescriptionlabel{#1}
}
\def\th@plain{
 \thm@notefont{}
 \itshape
}
\def\th@definition{
 \thm@notefont{}
 \normalfont
}

\g@addto@macro\th@definition{\thm@headpunct{}}
\g@addto@macro\th@plain{\thm@headpunct{}}
\makeatother

\usepackage{subfig}

\usepackage[final]{showkeys}

\usepackage{etoolbox}

\usepackage{mathrsfs}

\usepackage[titletoc]{appendix}
\usepackage[doc]{optional}

\usepackage{soul}

\usepackage{colortbl,booktabs,
multirow}
\usepackage{xcolor}

\usepackage{cancel}

\usepackage{empheq}

\definecolor{myblue}{rgb}{.8, .8, 1}
\newcommand*\mybluebox[1]{
\colorbox{myblue}{\hspace{1em}#1\hspace{1em}}}

\usepackage{hyperref}
\hypersetup{
colorlinks=true,
linkcolor=blue,
citecolor=blue,
filecolor=magenta,
urlcolor=cyan
}

\usepackage[
open,
openlevel=2,
atend,
numbered
]{bookmark}

\usepackage[capitalize, nameinlink, noabbrev]{cleveref}
\crefname{equation}{}{}
\crefname{chapter}{Chapter}{Chapters}
\crefname{item}{item}{items}
\crefname{figure}{Figure}{Figures}
\crefname{theorem}{Theorem}{Theorems}
\crefname{lemma}{Lemma}{Lemmas}
\crefname{proposition}{Proposition}{Propositions}
\crefname{corollary}{Corollary}{Corollarys}
\crefname{definition}{Definition}{Definitions}
\crefname{fact}{Fact}{Facts}
\crefname{example}{Example}{Examples}
\crefname{algorithm}{Algorithm}{Algorithms}
\crefname{remark}{Remark}{Remarks}
\crefname{note}{Note}{Notes}
\crefname{notation}{Notation}{Notations}
\crefname{case}{Case}{Cases}
\crefname{exercise}{Exercise}{Exercises}
\crefname{question}{Question}{Questions}
\crefname{claim}{Claim}{Claims}
\crefname{enumi}{}{}

\usepackage{pgf}

\usepackage{array}
\usepackage{tabu}

\allowdisplaybreaks

\numberwithin{equation}{section}

\renewcommand\theenumi{(\roman{enumi})}
\renewcommand\theenumii{(\alph{enumii})}

\spnewtheorem*{Proof}{Proof.}{\bf}{\rm}

\begin{document}

\title{Subtransversality and Strong CHIP of  Closed Sets in Asplund Spaces\thanks{Research  of the first author was supported by the Natural Science Foundation of Hebei Province (A2024201015), the Excellent Youth Research Innovation Team of Hebei University (QNTD202414) and the Innovation Capacity Enhancement Program-Science and Technology Platform Project, Hebei Province (22567623H). Research of the second author   was partially supported  by  FASIC and a public grant as part of the Investissement d'avenir project, reference ANR-11-LABX-0056-LMH, LabEx LMH.}}

\titlerunning{Subtransversality and Strong CHIP of Closed Sets in Asplund Spaces}

\author{Zhou Wei  \and Michel Th\'era \and Jen-Chih Yao}

\institute{Zhou Wei\at Hebei Key Laboratory of Machine Learning and Computational Intelligence \& College of Mathematics and Information Science, Hebei University, Baoding, 071002, China\\ \email{weizhou@hbu.edu.cn}\\
Michel Th\'era \at XLIM UMR-CNRS 7252, Universit\'e de Limoges, Limoges, France\\\email{michel.thera@unilim.fr}\\
Jen-Chih Yao \at Research Center for Interneural Computing, China Medical University Hospital,
China Medical University, Taichung, Taiwan \\ and 
\at Academy of Romanian Scientists, 50044 Bucharest, Romania\\ \email{yaojc@mail.cmu.edu.tw}
}

\date{Received: date / Accepted: date}

\maketitle

\begin{abstract}
In this paper, we mainly study subtransversality and two types of strong CHIP (given via Fr\'echet and limiting normal cones) for a collection of finitely many closed sets. We first prove characterizations of Asplund spaces in terms of subtransversality and intersection formulae of Fr\'echet normal cones. Several necessary conditions for subtransversality of closed sets are obtained via Fr\'echet/limiting normal cones in  Asplund spaces. Then, we consider subtransversality for some special closed sets in convex-composite optimization. In this frame we prove an  equivalence result on subtransversality, strong Fr\'echet CHIP and property (G) so as to extend a duality  characterization of subtransversality of finitely many closed convex sets via strong CHIP and property (G) to the possibly non-convex case.  As applications, we use these results on subtransversality and strong CHIP to study error bounds of  inequality systems  and give several dual criteria for error bounds via Fr\'echet normal cones and subdifferentials.




\end{abstract}

\section{Introduction and preliminary material}

\subsection{\textbf{ Definitions, notations and preliminary material}}
 
Let us begin with  some notations and definitions. Throughout the text, the notation employed   is fairly standard.  Let $\mathbb{X}$  be  a  real Banach space  and  $\mathbb{X}^\ast $  the dual space of continuous linear functionals on   $\mathbb{X}$.  We  use the notation $\langle\cdot,\cdot\rangle$ to indicate  the duality pairing between $\mathbb{X}$ and $\mathbb{X}^\ast$.
 As usual,  $\mathbf{B}_{\mathbb{X}}:= \{x\in \mathbb{X}:\Vert x\Vert \leq 1\}$  stands for the closed unit ball in $\mathbb{X}$ and for $x\in \mathbb{X},$ $r>0,$ $\mathbf{B}(x,r):= \{y\in \mathbb{X} : \Vert x-y\Vert <r\}$ refers to the open ball centered at $x$ with radius $r$, respectively. If $x\in \mathbb{X}$, we define the distance from $x$ to a subset 
$C$ of $\mathbb{X}$   by $$\mathbf{d}(x,C): = 
 \inf_{y\in C}\Vert x-y\Vert. $$
 By a cone we mean a nonempty set in $\mathbb{X}$ which is closed under nonnegative scalar multiplication. Given 
 a nonempty set $S\subseteq \mathbb{X}$,  the notation  $\cone S$  is used for the conical hull of $S$, i.e., the cone generated by $S$. The set $\cone S$     is  convex whenever $S$ is convex. 
 Suppose $S$  is a nonempty set in $\mathbb{X}$. Then the negative polar cone  of $S,$ written $ S^0$  is defined by $\{x^\ast \in \mathbb{X^\ast} : \langle x^\ast  ,x\rangle \leq 0\quad \text{for every } x\in S\}$.
{The symbol $\stackrel{w^\ast }\longrightarrow $ stands for   the weak$^\ast $ convergence and the symbol  $y \stackrel{A}\longrightarrow x $ means that $y$ tends to $x$  while staying in a given subset $A$ of  $\mathbb{X}$.
Let $A$ be a closed subset of $\mathbb{X}$ and $x\in A$. For  $\varepsilon\geq 0$, we denote by
$$
\hat N_{\varepsilon}(A, x):=\left\{x^\ast \in
X^\ast :\limsup\limits_{y\stackrel{A}\rightarrow x}\frac{\langle x^\ast ,
y-x\rangle}{\|y-x\|}\leq \varepsilon\right\}
$$
the set of $\varepsilon$-\textit{normals} to $A$ at $x$. When $\varepsilon=0$,  we will write  $\hat N(A,x)$ instead of $\hat N_{0}(A,x)$. The set $\hat N(A,x)$  is a closed convex cone  called the \textit{Fr\'echet normal cone to $A$ at $x$. } 

Given a multifunction $\Phi: \mathbb{X}\rightrightarrows \mathbb{X}^\ast $, the symbol
\begin{equation*}
\begin{array}r
\mathop{\rm Limsup}\limits_{y\rightarrow x}\Phi(x):=\Big\{x^\ast \in \mathbb{X}^\ast : \exists \ {\rm sequences} \ x_n\rightarrow x \ {\rm and} \ x_n^\ast \stackrel{w^\ast }\longrightarrow x^\ast  \ {\rm with}\  \\
x_n^\ast \in \Phi(x_n) \ {\rm for\ all \ } n\in \mathbb{N} \Big\}
\end{array}
\end{equation*}
signifies the  sequential Painlev\'{e}-Kuratowski outer/upper limit of $\Phi(x)$ as $y\rightarrow x$.
The  \textit{ limiting (Mordukhovich) normal cone}  of $A$ at $x$   is defined by
$$
N(A, x):=\mathop{\rm Limsup}_{y\stackrel A\rightarrow x,\, \varepsilon\downarrow
0}\hat N_{\varepsilon}(A, y).
$$
Thus $N(A,x)$ is the set of those $x^\ast \in \mathbb{X}^\ast $   such that 
there exists a sequence $\{(x_n,\varepsilon_n,x_n^\ast )\}$ in $A\times \mathbb{R}_+\times \mathbb{X}^\ast $ such that $(x_n,\varepsilon_n)\rightarrow(x,0)$, $x_n^\ast \stackrel{w^\ast }{\rightarrow}x^\ast $ and $x^\ast _n\in\hat{N}_{\varepsilon_n}(A,x_n)$ for each $n$. It is known from 
\cite{Mordukhovich} that
$\hat N(A,x)\subseteq N(A, x)$  with an equality when $A$ is convex. In this case these cones coincide  with  the normal cone (in the sense of the classical convex analysis):
$$
\hat{N}(A,a)=N(A, a)=\{x^\ast \in \mathbb{X}^\ast :\;\langle
x^\ast ,x-a\rangle\leq 0\;\; {\rm for\ all} \ x\in A\}.
$$
%
Recall that a Banach space is  Asplund if every continuous convex function on $\mathbb{X}$ is Fr\'echet differentiable at each point of a dense subset of $\mathbb{X}$ (see \cite{Phelps, Yost93} for definitions and their equivalences).  
An  important fact should be underlined:
any Banach space, with an equivalent norm which is Fr\'echet differentiable  off  zero, is Asplund. This property implies that every reflexive Banach space is Asplund since any such space admits an equivalent Fr\'echet differentiable norm off zero   (see \cite[Definition 2.2.3]{thibault}).
}
\vskip 2mm
Another important property of Asplund spaces is the   sequential representation of  the limiting normal cone established  by Mordukhovich and Shao \cite{MS}: 
\begin{equation*}
N(A, a)=\mathop{\rm Limsup}_{y\stackrel A\rightarrow
x}\hat N(A, y).
\end{equation*}
Equivalently,   $x^\ast \in N(A, x)$ if and only if there exist sequences $\{x_n\}$  and $\{x_n^\ast \}$ such that  $x_n\stackrel{A}\rightarrow x$, $x^\ast _n\stackrel{w^\ast }\rightarrow x^\ast $ and $x_n^\ast \in \hat N(A, x_n)$ for all $n$.
\vskip 2mm
Let  $f:\mathbb{X}\rightarrow \mathbb{R}\cup\{+\infty\}$ be a proper lower semicontinuous  extended real-valued function. We denote as usual by
$$
\dom f:=\{y\in \mathbb{X}:f(y)<+\infty\} \ \ {\rm and} \ \ \epi f:=\{(x,\alpha)\in\mathbb{X}\times\mathbb{R}: f(x)\leq \alpha\}
$$ 
the \textit{domain } and the\textit{ epigraph} of $f$, respectively. The  \textit{Fr\'{e}chet} and \textit{ limiting} subdifferentials of $f$ at $x\in\dom f$ are the subsets of $\mathbb{X}^\ast $  which are  respectively  defined by
$$
\hat{\partial}f(x):=\{x^\ast \in \mathbb{X}^\ast : (x^\ast ,-1)\in \hat N(\epi f, (x, f(x)))\} 
$$
and
$$
 \partial f(x):=\{x^\ast \in \mathbb{X}^\ast : (x^\ast ,-1)\in N(\epi f, (x, f(x)))\}.
$$
Equivalently, 
$$
\hat{\partial}f(x)=\left\{x^\ast \in \mathbb{X}^\ast : \liminf_{y\rightarrow x}\frac{f(y)-f(x)-\langle x^\ast , y-x\rangle}{\|y-x\|}\geq 0\right\},
$$
and  when $\mathbb{X}$ is  Asplund,  Mordukhovich and Shao \cite{MS} have proved that
\begin{equation*}
\partial f(x)=\mathop{\rm Limsup}_{y\stackrel f\rightarrow
x}\hat \partial f(y),
\end{equation*}
where $y\stackrel f\rightarrow x$ means that $y\rightarrow x$ and $f(y)\rightarrow f(x)$. Thus $x^\ast \in\partial f(x)$ if and only if there exist $x_n\stackrel{f}\rightarrow x$ and $x^\ast _n\stackrel{w^\ast }\rightarrow x^\ast $ such that $x_n^\ast \in \hat \partial f(x_n)$ for all $n$.

\subsection{\textbf{ A brief survey on the  notions of  convex feasibility problem, strong conical hull intersection property
(strong CHIP) and property (G),  subtransversality.}}

The concepts of \textit{subtransversality, strong conical hull intersection property
(strong CHIP) and property (G) } for a collection of closed convex sets have   played important roles  in various branches of approximation theory and optimization. They are used  in the  convex feasibility problem, constrained approximation, Fenchel duality, systems of convex inequalities and error bounds (see  \cite{BDL2005,BBL,BBT,HW,LJ,LN,LN1,LN2,LNP,WH} and references therein).  
%

Let us briefly survey these concepts.
\vskip 2mm
\textbf{Convex feasibility problem (CFP)}
\vskip 2mm
The  \textit{convex feasibility problem} is the problem  of finding a common point in   a collection of finitely many closed convex intersecting sets (when this intersection is nonempty).    This  general formulation is of  practical  interest   in diverse areas of mathematics and physical sciences   and  is therefore  central in various applications such as for instance, 
 solution of convex inequalities, minimization of convex non-smooth functions, medical imaging, computerized tomography and electron microscopy. CFP   often solved iteratively by projection algorithms,  has been the subject  of  extensive research to ensure a better norm convergence of the sequences produced by  projection-based algorithms (cf. \cite{BB,BB1,BBL97,KL1997}). We refer  the reader  to the   survey of Bauschke \& Borwein \cite{BB1} or to   the contribution of  Combettes \cite{Combettes96}  and the references therein. In particular Combettes \cite{Combettes96}  gives a concrete application  to  the image recovery problem.  
CFP has also a connection with the theory of error bounds. 
\vskip 2mm
\textbf{Conical Hull Intersection Property (CHIP)}
\vskip 2mm
Chui, Deutsch and Ward \cite{CDW1,CDW2} introduced the concept of  \textit{conical hull intersection property } (CHIP) when studying the problem of constrained interpolation from a convex subset or cone.  Given a finite family of closed convex sets $\{C_1,\cdots, C_m\}$  of $\mathbb{X}$ such that $C:=\bigcap_{i=1}^mC_i\neq\emptyset$,   we say that the CHIP holds at $x\in C$  if
\begin{empheq}[box =\mybluebox]{equation*}T(C, x) = \bigcap_{i=1}^mT(C_i, x),\end{empheq}
where $T(C, x)$, the closure of $\cone (C-x)$, stands for the \textit{tangent cone} to $C$ at $x$. Note that $T(C, x)= N(C,x)^{\circ}:=\{h\in\mathbb{X}: \langle x^*,h\rangle\leq 0,\ \forall x^*\in N(C,x)\}$.
\vskip 2mm
\textbf{Strong Conical Hull Intersection Property (Strong CHIP)}
\vskip 2mm

It is well-known that important methods to solve the CFP in Hilbert spaces are based on some projection-based algorithms (including serial 	methods, simultaneous methods, extrapolated relaxation method, and subgradient 	methods) and usually the sequences obtained from these algorithms have the weak-convergence property, not strong-convergence in general. For ensuring a better 	convergence property, Deutsch, Li and Ward \cite{DeuLiWar99} introduced strong CHIP, a stronger version of CHIP, to study an unconstrained reformulation of the constrained best approximation. Given a finite family of closed convex sets $\{C_1,\cdots, C_m\}$  of $\mathbb{X}$ such that $C:=\bigcap_{i=1}^mC_i\neq\emptyset$,   we say that the strong CHIP holds at $x\in C$  if
  \begin{empheq}[box =\mybluebox]{equation*}N(C, x) = \sum_{i=1}^mN(C_i, x).\end{empheq}
It is known that the pair $(C,D)$  of closed and convex subsets of a Banach space $\mathbb{X} $ has the strong CHIP if and only if the sum rule formula applies for the indicator functions $i_C$ and $i_D$ of $C$ and $D$, where $i_C (x)=0$  if $x\in C$ and $+\infty$ otherwise.  
  According to   a classical result due to  Moreau \cite{moreau}, the strong CHIP is also equivalent to saying  that the infimal convolution of their support functionals  is exact.
 Later, it was  proved by  Deutsch \cite{D} that strong CHIP is a geometric version of the basic constraint qualification for constrained optimization problems. Deutsch, Li and Ward \cite{DeuLiWar99} further proved that strong CHIP is necessary and sufficient for some ``perturbation property". Deutsch, Li and Swetits \cite{DLS} showed that strong CHIP is an essential property when establishing duality relations between certain convex optimization problems. Jeyakumar further \cite{J1} proved that  strong CHIP is equivalent to a complete Lagrange multiplier characterization of optimality for convex programming model problems. Ernst and Th\'era \cite{ET2007} considered strong CHIP for the pair of closed and convex sets, and used the strong CHIP to prove that the validity of the converse of Moreau's theorem is associated with the absence of half-strips within the boundary of the  concerned convex sets.
Let us point out   an important fact:  if  strong CHIP holds, then  when minimizing a convex function  $f$ over the intersection of two convex sets $C$ and $D$,   the first order optimality condition  becomes 
\begin{empheq}[box =\mybluebox]{equation*}\
0\in \partial f(x) + N(C, x)+ N(D,x).\end{empheq}
Deutsch proved in \cite{Deutsch98} that strong CHIP is the weakest constraint qualification under which a minimizer of a convex function$f : C\cap  D \to \R$ can be characterized using the subdifferential of $f$  at $x$ and the normal cones of $C$  and $D$  at $x$.

\vskip 2mm
\textbf{Property (G)}
\vskip 2mm
  In the early 1970s, Jameson \cite{J} introduced the so-called property (G). This property   provides  a sufficient condition for the sum of   two closed convex cones  $C $ and $D $   to be closed and  allows   to establish a duality theory for two closed convex cones. Given a finite family of closed convex cones $\{P_1,\cdots,P_m\}$ of $\mathbb{X}$, we say that $\{P_1,\cdots,P_m\}$ has property (G) if there exists $\tau\in (0,+\infty)$ such that
  
  	\begin{empheq}[box =\mybluebox]{equation*}
  		\Big(\sum_{i=1}^mP_i\Big)\cap \mathbf{B}_{\mathbb{X} }\subseteq \tau \sum_{i=1}^m(P_i\cap \mathbf{B}_{\mathbb{X} }).
  	\end{empheq}
 See \cite{BZ2005}.



\vskip 2mm
\textbf{Subtransversality}
\vskip 2mm

Given a collection of finitely many closed sets $\{A_1,\cdots,A_m\}$  in a Banach space $\mathbb{X}$, recall that the collection  $\{A_1,\cdots,A_m\}$ is said to be  \textit{  subtransversal } at $\bar x\in\bigcap_{i=1}^mA_i$ if there exists $\tau\in(0,+\infty)$ such that
\begin{empheq}[box =\mybluebox]{equation}\label{2}
	\mathbf{d}\Big(x,\bigcap_{i=1}^mA_i\Big)\leq \tau\sum_{i=1}^m\mathbf{d}(x, A_i)\ \ {\rm for\ all} \ x\ {\rm close\ to} \ \bar x.
\end{empheq}

\textit{Subtransversality} is a well-known notion in mathematical programming and approximation theory. It is known from \cite[7.1.3 Comments]{Ioffe2017} that inequality \eqref{2} was in reality introduced by Dolecki in a very different context (see\cite{Dolecki1982}). This property later appeared in \cite{Ioffe1989} as a qualification condition in the calculus of normal cones and was called a ``metric qualification condition" in \cite{IP1996};  it appears also as ``linear coherence"
in \cite[Theorem 4.7.5]{Penot2013}.  For convex sets, the inequality was also introduced in \cite{BB} under the name ``linear regularity". Subsequently, Bauschke and  Borwein \cite{BB1} used this notion to study the problem of finding the projection from a  given point to the  intersection of finitely many closed convex sets, and proved that subtransversality  plays a key role in establishing a linear convergence rate of the iterates generated by the cyclic projection algorithm.  We refer the reader to \cite{KLT2017,KLT2018,BCK2020} for a survey of the terminology and a comparison with similar or related concepts available in the literature. Since it is closely linked to many optimization  issues,  subtransversality has been extensively studied by many authors (cf. \cite{BDL2005,BB,BBT,LP,LNP,W,WZ} and references therein). 
\vskip 2mm

\vskip 2mm

\textbf{\textbf{Metric {Subregularity}}}
\vskip 2mm

Let $F:\mathbb{X}\rightrightarrows \mathbb{Y}$ be a  set-valued mapping between two Banach spaces $ \mathbb{X}$  and $\mathbb{Y}$ and $\bar x\in \mathbb{X}$. Recall that $F$ is said to be \textit{metrically subregular} at $\bar x$ for $\bar y\in F(\bar x)$, if there exists  a constant $\kappa>0$ such that  
\begin{empheq}[box =\mybluebox]{equation*}
 \mathbf{d}(x, F^{-1}(\bar y))\leq\kappa \mathbf{d}(\bar y, F(x))\ \  {\rm for\ all} \ x\ {\rm close\ to} \ \bar x.
\end{empheq}
Metric subregularity  has been thoroughly studied in various contexts in variational analysis and has been the object of a very active  area of research in both theory and applications (see \cite{DonRoc09,IO2008} and \cite{Kru15,Penot2013} for its relationship with the theory of error bounds).  

Using the $\ell^1$-norm on $\mathbb{X}^n$, let us consider the set-valued mapping $F:\mathbb{X}\rightrightarrows \mathbb{X}^n$  defined as follows:
$$F(x):= (A_1-x)\times (A_2-x)\times \cdots \times (A_n-x).$$ 
We easily observe that  
$$ \bigcap_{i=1}^mA_i  = F^{-1}(0) = \{x\in \mathbb{X}: 0\in F(x)\},$$  
$$  \mathbf{d}\Big(x,\bigcap_{i=1}^mA_i\Big) = \mathbf{d}(x, F^{-1}(0) )\quad  \text{and} \quad  \mathbf{d}(0, F(x)) = \sum_{i=1}^n \mathbf{d}(x,A_i).$$
Then  subtransversality  which satisfies  \eqref{2} is equivalent to  the  metric subregularity of $F$ at $(0,\bar x)$. It should be noted that this equivalence also holds even if the convexity assumption is dropped.  We refer  to \cite{NT2004} for details and results. Such equivalence was also given in \cite[Proposition 1]{KLT2017} and \cite[Theorem 1]{KLT2018}. 

  For closed convex cones, subtransversality is equivalent to property (G) (see \cite[Fact 2.11]{BBT}).
 In 1996, Bauschke and Borwein \cite{BB1} considered subtransversality when solving convex feasibility problems. Afterwards Pang \cite{P} and Lewis and Pang \cite{LP} established dual conditions for subtransversality of finitely many closed convex sets. In 1999, using convex analysis techniques, Bauschke, Borwein, and Li   \cite{BBL} discussed  the close relationship among subtransversality, strong CHIP and property (G) for a collection of finitely many closed convex sets in  Euclidean spaces. In 2005, Bakan, Deutsch and Li \cite{BDL2005} used various normal properties to study property (G), strong CHIP and subtransversality for a finite collection of convex sets in  Hilbert spaces. The  normal property, property (G), strong CHIP and subtransversality for an infinite collection of closed convex sets in a Banach space was further considered in \cite{LN2} and \cite{WH}. The following well-established theorem is a relevant result on subtransversality, strong CHIP and property (G) for a collection of closed convex sets (cf. \cite{BDL2005,NY,ZN2004,LN2,WH}):

\medskip

\noindent{\bf Theorem A.} {\it
 Let $\{A_1,\cdots, A_m\}$ be a collection of closed convex sets in a Banach space $\mathbb{X}$ such that  $A:=\bigcap_{i=1}^mA_i$ is nonempty. Then the following statements are equivalent:}
\begin{itemize}
  \item [\rm (i)] {\it $\{A_1,\cdots,A_m\}$ is subtransversal at $\bar x\in A$.}
  \item [\rm (ii)] {\it There exist $\eta,\delta\in (0,+\infty)$ such that for any $x\in A\cap \mathbf{B}(\bar x,\delta)$, $\{A_1,\cdots,A_m\}$ has the strong CHIP at $x$ and}
  $$
  \inf\left\{\sum_{i=1}^{m}\|x_i^\ast \|:x_i^\ast \in N(A_i,x) \ {\it and} \  x^\ast =\sum_{i=1}^{m}x_i^\ast \right\}\leq \eta\|x^\ast \|
  $$
  {\it holds for all $x^\ast \in N\left(A, x\right)$}.
  \item [\rm (iii)] {\it There exists $\tau,\delta\in (0,+\infty)$ such that for any $x\in A\cap \mathbf{B}(\bar x,\delta)$,  $\{A_1,\cdots,A_m\}$ has the strong CHIP at $x$ and $\{N(A_1,x),\cdots,N(A_n,x)\}$ has property (G) with constant $\tau$.}
\end{itemize}

 \medskip
 
This important theorem shows the equivalence relation among subtransversality, strong CHIP and property (G). A natural question that arises is whether  similar results still hold  for general closed sets. More precisely,  whether subtransversality could also be characterized by  strong CHIP and property (G) in the absence of convexity assumption?
Unfortunately, the answer to this question is negative once the convexity assumption is dropped (see \cref{re3.3}(b) below). The first difficulty to deal with the question under the absence of convexity  is to choose one suitable kind of normal cones and then consider the corresponding concepts of strong CHIP and property (G) in terms of such normal cones. Along this line, 
   our goal is to  deepen the relation between  subtransversality,   strong CHIP and   property (G) for a collection of a finite number of closed sets in Asplund spaces.
\subsection{\textbf{Organization of the paper}}
This paper is organized as follows. Several preliminaries  and known
results will be given in Section 2. In Section 3, we first prove that a   Banach space is   Asplund if and only if the intersection formulae of Fr\'echet normal cones for every nonempty closed sets is  implied by  subtransversality. Several necessary conditions for subtransversality of closed sets are given via Fr\'echet/limiting normal cones. 
Then we study subtransversality for some special closed sets in convex-composite optimization and prove that  subtransversality  is equivalent to  the strong Fr\'echet CHIP and property (G)  when $\mathbb{X} $ is Asplund (see \cref{th3.4}). 
In Section 4 we apply the results obtained in Section 3 to the study of error bounds of inequality systems. Conclusions of this paper are presented in Section 5.
\section{Preliminaries}
The following lemma  concerns the fuzzy sum rule for Fr\'echet subdifferentials in Asplund spaces which  was proved by Fabian \cite[Theorem 3]{F}. Readers could also consult \cite[Proposition 2.7]{MS} and \cite[Theorem 2.33]{Mordukhovich} for details.

\begin{lemma}[fuzzy sum rule for Fr\'echet subdifferentials in Asplund spaces]\label{lem2.1}
Let $\mathbb{X}$ be an Asplund space and $f_i:\mathbb{X}\rightarrow \mathbb{R}\cup\{+\infty\},i=1,2,$ be proper lower semicontinuous functions one of which is locally Lipschitz at $\bar x\in\dom f_1\cap \dom f_2$. Then for any $\varepsilon>0$, one has
\begin{eqnarray*}
\hat \partial(f_1+f_2)(\bar x)\subset\bigcup\{\hat \partial f_1(x_1)+\hat \partial f_2(x_2):x_i\in \mathbf{B}(\bar x,\varepsilon),|f_i(x_i)-f(\bar x)|<\varepsilon,
 i=1,2\}+\varepsilon \mathbf{B}_{\mathbb{X}^\ast }.
\end{eqnarray*}
\end{lemma}

The next  lemma is taken from \cite[Proposition 2.1]{NT2001} which will be   used in the sequel.

\begin{lemma}\label{lem2.2a}
 Let $f:\mathbb{X}\rightarrow \mathbb{R}\cup\{+\infty\}$ be a proper lower semicontinuous function and $(\bar x,\alpha)\in \epi f$. Then for any $\lambda\not =0$, the following equivalence is true in any Banach space $\mathbb{X}$:
$$
(x^\ast ,-\lambda)\in\hat N(\epi f, (\bar x,\alpha))\ \Longleftrightarrow\ \lambda>0, \alpha=f(\bar x), \frac{x^\ast }{\lambda}\in\hat\partial f(\bar x).
$$
\end{lemma}

Lemma \ref{lem2.2} (see \cite[Lemma 3.6]{NT2001}  and also also \cite[Lemma 3.7]{AusDanThi05})  and Lemma \ref{lem2.3}  (see  \cite[Theorem 3.1]{ZN} ) are  useful in our analysis.

\begin{lemma}\label{lem2.2}
Let $\mathbb{X}$ be an Asplund space and $A$ be a nonempty closed subset of
$\mathbb{X}$. Let $x\notin A$ and $x^\ast \in \hat{\partial}\mathbf{d}(\cdot,
A)(x)$. Then, for any $\varepsilon>0$ there exist $a\in A$ and
$a^\ast \in \hat N(A, a)$ such that
$$
\|x-a\|<\mathbf{d}(x,
A)+\varepsilon\,\,\,\,{\it and}\,\,\,\,\|x^\ast -a^\ast \|<\varepsilon.
$$
\end{lemma}


\begin{lemma}\label{lem2.3}
Let $\mathbb{X}$ be  an Asplund space,  $A$ be a nonempty closed
subset of $\mathbb{X}$ and $x\not\in A$. Then for any $\beta\in (0,\;1)$
there exist $z\in A$ and $z^\ast \in
\hat{N}(A, z)$ with $\|z^\ast \|=1$ such that
$$\beta\|x-z\|<\min\{\mathbf{d}(x, A),\;\langle z^\ast ,x-z\rangle\}.$$
\end{lemma}

In what follows, given  a collection of finitely many closed sets, we consider the notions of  strong Fr\'echet and limiting CHIP, respectively.

\begin{definition}\label{def2.1}

Let $\{A_1,\cdots,A_m\}$ be a collection of finitely many closed sets in $\mathbb{X}$ such that $A:=\bigcap_{i=1}^mA_i$ is nonempty, and let $\{P_1,\cdots,P_m\}$ be a collection of finitely many closed cones in $\mathbb{X}^\ast $. Recall that
\begin{itemize}
\item[\rm (i)] $\{A_1,\cdots,A_m\}$ is said to have the \textit{strong Fr\'echet conical hull intersection property}  (\textbf{strong Fr\'echet CHIP)}  at $x\in A$ if
\begin{empheq}[box =\mybluebox]{equation}\label{2.2}
\hat N(A,x)\subseteq \sum_{i=1}^m\hat N(A_i,x)
\end{empheq}
and $\{A_1,\cdots,A_m\}$ is said to have \textit{the strong Fr\'echet CHIP}  if $\{A_1,\cdots,A_m\}$ has the strong Fr\'echet CHIP at all points of $A$;

\item[\rm (ii)] $\{A_1,\cdots,A_m\}$ is said to   have  the \textbf{strong limiting CHIP} at $x\in A$ if
\begin{empheq}[box =\mybluebox]{equation*}\label{2.3}
N(A,x)\subseteq \sum_{i=1}^m N(A_i,x)
\end{empheq}
and $\{A_1,\cdots,A_m\}$ is said to have the strong limiting CHIP if $\{A_1,\cdots,A_m\}$ has the strong limiting CHIP at all points of $A$.
\item[\rm (iii)] $\{P_1,\cdots,P_m\}$ is said to have the \textbf{property ($G$)} if there exists $\tau\in (0,+\infty)$ such that
\begin{empheq}[box =\mybluebox]{equation*}
  \Big(\sum_{i=1}^mP_i\Big)\cap \mathbf{B}_{\mathbb{X}^\ast }\subseteq \tau \sum_{i=1}^m(P_i\cap \mathbf{B}_{\mathbb{X}^\ast }).
\end{empheq}
\end{itemize}
\end{definition}

\begin{remark}\label{rem2.1}
\begin{itemize}
\item[\rm (i)] Since $\hat N(A,x)$ is always a convex cone and contains $\hat N(A_i,x)$ for each $i$, the requirement \eqref{2.2} for the strong Fr\'echet CHIP is equivalent to
\begin{empheq}[box =\mybluebox]{equation*}\label{2.5}
\hat N(A,x)= \sum_{i=1}^m \hat N(A_i,x).
\end{empheq}
\item[\rm (ii)] We note that   the strong Fr\'echet/limiting CHIP refers to calculus of Fr\'echet/limiting normal cones (cf. \cite[Lemma 3.1 and Theorem 3.4]{Mordukhovich}). For the case when $\{A_1,\cdots,A_m\}$ is a collection of finitely many closed convex sets, the notions of strong Fr\'echet CHIP and strong limiting CHIP coincide and reduce to the strong CHIP studied in the convex optimization.
\end{itemize}
\end{remark}

We conclude this section with the following example to show that strong Fr\'echet CHIP and strong limiting CHIP for non-convex closed sets could be valid.

\begin{example}	\begin{itemize}
\item [\rm(a)] 	Let $\mathbb{X}:=\mathbb{R}^2$, $\bar x:=(0,0)\in \mathbb{X}$, $A_1:=\mathbb{R}_-\times\mathbb{R}_-$ and $A_2:={\rm cl}(\mathbb{X}\backslash A_1)$. Then $$A_1\cap A_2=(\mathbb{R}_-\times\{0\})\cup(\{0\}\times\mathbb{R}_-).$$ 
	Then $\hat N(A_1\cap A_2,\bar x)=\hat N(A_1,\bar x)=\mathbb{R}_+\times\mathbb{R}_+$, $\hat N(A_2,\bar x)=\{(0,0)\}$ and thus $\{A_1,A_2\}$ has the strong Fr\'echet CHIP at $\bar x$.
	\item [\rm(b)] Let $\mathbb{X}:=\mathbb{R}^2$, $\bar x:=(0,0)\in \mathbb{X}$, $A_1:=(\mathbb{R}_-\times\mathbb{R}_+)\cup\{(s,t)\in \mathbb{R}^2: s+t=0\}$ and
	$$
	A_2:=\left\{(s, t)\in\mathbb{R}^2_+: (s-1)^2+t^2\leq 1\ \ {\rm and} \ \ s^2+(t-1)^2\leq 1\right\}.
	$$
	Then $A_1\cap A_2=\{\bar x\}$ and by computation, one has $\hat N(A_1\cap A_2,\bar x)=N(A_1\cap A_2,\bar x)=\mathbb{X}$,
	$$
	\hat N(A_2,\bar x)=N(A_2,\bar x)=\mathbb{R}_-\times\mathbb{R}_- \ \ {\rm and} \ \ \hat N(A_1,\bar x)=\{(0,0)\}
	,
	$$
	which means that strong Fr\'echet CHIP of $\{A_1,A_2\}$ at $\bar x$ is violated. However,
	$$
	 N(A_1,\bar x) = ([0, +\infty)\times\{0\})\cup(\{0\}\times (-\infty, 0])\cup\{(t,t):t\in\mathbb{R}\}.
	$$
	This implies that  $\{A_1,A_2\}$ has the strong limiting CHIP at $\bar x$.
\end{itemize}
\end{example}

\setcounter{equation}{0}

\section{Subtransversality and Strong CHIP of Closed Sets in \tcv{an} Asplund Space}

In this section, we mainly study subtransversality for a collection of finitely many closed sets in  Asplund spaces. We aim to provide necessary and/or sufficient conditions for subtransversality via strong Fr\'echet/limiting CHIP and property (G). We begin with the definition of subtransversality.

Let $\{A_1,\cdots,A_m\}$ be a collection of finitely many closed sets in $\mathbb{X}$ such that the intersection $A:=\bigcap_{i=1}^mA_i$ is nonempty. Recall that the collection $\{A_1,\cdots,A_m\}$ is said to be  subtransversal at $\bar x\in A$ if there exist $\tau, r\in(0,+\infty)$ such that
\begin{empheq}[box =\mybluebox]{equation}\label{3.3}
  \mathbf{d}(x, A)\leq \tau\sum_{i=1}^m\mathbf{d}(x, A_i)\ \ \forall x\in \mathbf{B}(\bar x, r).
\end{empheq}

We first prove the following result on necessary conditions for subtransversality of finitely many closed sets in an Asplund space.

\begin{theorem}\label{th3.2}
	Let $\{A_1,\cdots,A_m\}$ be a collection of finitely many closed sets in an Asplund space $\mathbb{X}$ such that $A:=\bigcap_{i=1}^mA_i$ is nonempty. Suppose that $\{A_1,\cdots,A_m\}$ is  subtransversal at $\bar x\in A$. Then there exist $\tau,\delta\in(0, +\infty)$ such that for any $\varepsilon>0$, one has
	\begin{equation}\label{3.11}
		\hat N(A,x)\cap \mathbf{B}_{\mathbb{X}^\ast }\subset\bigcup\Big\{\tau\sum_{i=1}^m\hat N(A_i, x_i)\cap(1+\varepsilon)\mathbf{B}_{\mathbb{X}^\ast }: x_i\in A_i\cap \mathbf{B}(x, \varepsilon), i=1,\cdots,m \Big\}+\varepsilon \mathbf{B}_{\mathbb{X}^\ast }
	\end{equation}
	holds for all $x\in \mathbf{B}(\bar x,\delta)\cap A$.
\end{theorem}

{\bf Proof.} Since $\{A_1,\cdots,A_m\}$ is  subtransversal at $\bar x$, then there exist $\tau, r\in(0,+\infty)$ such that \eqref{3.3} holds. Let $\delta:=\frac{r}{2}$ and $\varepsilon>0$. Take any $x\in \mathbf{B}(\bar x,\delta)\cap A$ and $x^\ast \in\hat N(A,x)\cap \mathbf{B}_{\mathbb{X}^\ast }$. Choose $\varepsilon_1>0$ such that $(m+4+\tau m)\varepsilon_1<\varepsilon$. Note that $
x^\ast \in \hat N(A,x)\cap \mathbf{B}_{\mathbb{X}^\ast }=\hat \partial \mathbf{d}(\cdot, A)(x)
$ and thus there exists $\delta_1>0$ such that
$$
\langle x^\ast ,u-x\rangle\leq \mathbf{d}(u, A)+\varepsilon_1\|u-x\|\ \ \forall u\in \mathbf{B}(x, \delta_1).
$$
By \eqref{3.3}, one has
$$
\langle {x}^\ast ,u-x\rangle\leq \tau\sum_{i=1}^m\mathbf{d}(u, A_i)+\varepsilon_1\|u-x\|\ \ \forall u\in \mathbf{B}(x, \delta_1).
$$
This implies that

$$
x^\ast \in\hat\partial \Big(\tau\sum_{i=1}^m\mathbf{d}(\cdot,A_i)+\varepsilon_1\|\cdot-x\|\Big)(x).
$$
By virtue of \cref{lem2.1,lem2.2}, there exist $u_i\in \mathbf{B}(x,\varepsilon_1)$, $u_i^\ast \in\hat\partial \mathbf{d}(\cdot,A_i)(u_i)$, $x_i\in A_i$ and $x^\ast _i\in \hat N(A_i,x_i) (1\leq i\leq m)$ such that
\begin{equation}\label{3.12}
	\|x_i-u_i\|<\mathbf{d}(u_i, A_i)+\varepsilon_1, \|x_i^\ast -u_i^\ast \|<\varepsilon_1
\end{equation}
and
\begin{equation}\label{3.13}
	x^\ast \in\tau\sum_{i=1}^mu_i^\ast +\varepsilon_1\mathbf{B}_{\mathbb{X}^\ast }+m\varepsilon_1\mathbf{B}_{\mathbb{X}^\ast }\subset \tau\sum_{i=1}^mx_i^\ast +(m+1+\tau m)\varepsilon_1\mathbf{B}_{\mathbb{X}^\ast }.
\end{equation}
For any $i$, one has $\|u_i^\ast \|\leq 1$ and thus $\|x_i^\ast \|\leq 1+\varepsilon_1<1+\varepsilon$ by \eqref{3.12}. This and \eqref{3.13} imply that
$$
x^\ast \in\big(\tau\sum_{i=1}^m\hat N(A_i, x_i)\cap(1+\varepsilon)\mathbf{B}_{\mathbb{X}^\ast }\big)+\varepsilon \mathbf{B}_{\mathbb{X}^\ast }.
$$
(thanks to the choice of $\varepsilon_1$). Since $\|x_i-x\|\leq \|x_i-u_i\|+\|u_i-x\|<3\varepsilon_1<\varepsilon$, it follows that \eqref{3.11} holds. The proof is complete. \hfill$\Box$

\medskip

It is known from \cite[Theorem 3.1 and Corollary 3.2]{NT2001} that necessary or sufficient conditions via  the Fr\'echet subdifferential of $\sum_{i=1}^m\mathbf{d}(\cdot, A_1)$ are provided to ensure the subtransversality (under the name of {\it metric inequality} therein). Further,  Ngai and Th\'era \cite[Theorem 3.8]{NT2001} gave characterizations of Asplund spaces in terms of conditions ensuring the subtransversality and the intersection formulae. In particular, the result is presented as follows:

\medskip

\noindent{\bf Theorem B.} \textit{ Let $\mathbb{X}$ be a Banach space. Then the following statements are equivalent:
\begin{itemize}
  \item [\rm{(i)}] $\mathbb{X}$ is an Asplund space;
  \item [\rm{(ii)}] for every nonempty closed subsets $A_1,\cdots, A_m$ of $\mathbb{X}$,  all but one of which are sequentially normally compact at $\bar x\in A_1\cap\cdots\cap A_m$ and satisfying the following condition:
      $$
      x_i^\ast \in N(A_i,\bar x), i=1,\cdots, m \ \ \&  \ \ x_1^\ast +\cdots+x_m^\ast =0 \Rightarrow x_1^\ast =\cdots=x_m^\ast =0,
      $$
      then $\{A_1,\cdots,A_m\}$ is  subtransversal at $\bar x$;
  \item [\rm{(iii)}] for every nonempty closed subsets $A_1,\cdots, A_m$ of $\mathbb{X}$ being  subtransversal at $\bar x\in A_1\cap\cdots\cap A_m$,  $\{A_1,\cdots,A_m\}$ has the strong limiting CHIP at $\bar x$.
\end{itemize}
}

\medskip

Similar to Theorem B, we prove the following theorem on characterizations of Asplund spaces given by the subtransversality and intersection formulae of Fr\'echet normal cones. This theorem also provides necessary conditions in terms of Fr\'echet normal cones for the subtransversality of finitely many closed sets in an  Asplund space.

\medskip

\begin{theorem}\label{th3.1}
Let $\mathbb{X}$ be a Banach space. Then the following statements are equivalent:
\begin{itemize}
  \item [\rm{(i)}] $\mathbb{X}$ is an Asplund space;
  \item [\rm{(ii)}] for every family $\{A_1,\cdots, A_m\}$ of nonempty closed  subsets of $\mathbb{X}$  which is   subtransversal at $\bar x\in A_1\cap\cdots\cap A_m$, there exists $\delta\in (0,+\infty)$ such that for any $\varepsilon>0$,  the inclusion 
        \begin{equation}\label{3.2}
  \hat N(A_1\cap\cdots\cap A_m, x)\subset\bigcup\Big\{\sum_{i=1}^m\hat N(A_i, x_i): x_i\in A_i\cap \mathbf{B}(x, \varepsilon), i=1,\cdots,m\Big\}+\varepsilon \mathbf{B}_{\mathbb{X}^*}
\end{equation}
      holds for all $x\in A_1\cap\cdots\cap A_m\cap \mathbf{B}(\bar x, \delta)$;
  \item [\rm{(iii)}]  for every family $\{A_1,\cdots, A_m\}$  of nonempty closed  subsets of $\mathbb{X}$ which is   subtransversal at $\bar x\in A_1\cap\cdots\cap A_m$
     and for any $\varepsilon>0$, one has  that \eqref{3.2} holds for $x=\bar x$. 
\end{itemize}
\end{theorem}

{\bf Proof.} (i)$\Rightarrow$(ii): Let $\{A_1,\cdots,A_m\}$ be  a collection of closed sets in $\mathbb{X}$ and denote $A:=\bigcap_{i=1}^mA_i$. Suppose that   the collection  $\{A_1,\cdots,A_m\}$ is  subtransversal at $\bar x$. Then there exist $\tau, r\in(0,+\infty)$ such that \eqref{3.3} holds. Let $\delta:=\frac{r}{2}$ and $\varepsilon>0$.
Take any $x\in A\cap \mathbf{B}(\bar x, \delta)$.
 Let $x^\ast \in \hat N(A,x)\backslash\{0\}$. Choose $\varepsilon_1\in(0, \varepsilon)$ such that $\|x^\ast \|\varepsilon_1<\varepsilon$. By virtue of \cref{th3.2} as well as its proof, one has
 \begin{equation*}
 	\frac{x^\ast}{\|x^\ast\|}\in\bigcup\Big\{\tau\sum_{i=1}^m\hat N(A_i, x_i)\cap(1+\varepsilon_1)\mathbf{B}_{\mathbb{X}^\ast }: x_i\in A_i\cap \mathbf{B}(x, \varepsilon_1), i=1,\cdots,m \Big\}+\varepsilon_1 \mathbf{B}_{\mathbb{X}^\ast }.
 \end{equation*}
This and the choice of $\varepsilon_1$ imply that
  \begin{eqnarray*}
 x^*&\in&\bigcup\Big\{\tau \Vert x^*\Vert\sum_{i=1}^m\hat N(A_i, x_i)\cap(1+\varepsilon_1)\mathbf{B}_{\mathbb{X}^\ast }: x_i\in A_i\cap \mathbf{B}(x, \varepsilon_1), i=1,\cdots,m \Big\}+\|x^*\|\varepsilon_1 \mathbf{B}_{\mathbb{X}^\ast }\\
 	&\subseteq&\bigcup\Big\{\sum_{i=1}^m\hat N(A_i, x_i): x_i\in A_i\cap \mathbf{B}(x, \varepsilon), i=1,\cdots,m \Big\}+\varepsilon \mathbf{B}_{\mathbb{X}^\ast }
 \end{eqnarray*}
and consequently  \eqref{3.2} holds. 

Note that (ii)$\Rightarrow$(iii) follows immediately and it remains to prove (iii)$\Rightarrow$(i). 

Suppose on the contrary that $\mathbb{X}$ fails to be Asplund.  Following \cite{MW2000,FM1998}, we can represent $\mathbb{X}$ in the form $\mathbb{X}=\mathbb{Y}\times \mathbb{R}$ with the norm $\|(y,\alpha)\|:=\|y\|+|\alpha|$ for any $x=(y,\alpha)\in \mathbb{X}$. Then $\mathbb{Y}$ is not an Asplund space. By virtue of \cite[Theorem 1.5.3]{DGZ1993} (also \cite[Theorem 2.1]{FM1998}), there exists an equivalent norm $|||\cdot|||$ on $\mathbb{Y}$ and $\gamma>0$ such that 
\begin{equation}\label{3.05}
\frac{1}{2}\|y\|\leq|||y|||\leq \|y\|\ \ {\rm and} \ \  \limsup_{h\rightarrow 0}\frac{|||y+h|||+|||y-h|||-2|||y|||}{\|h\|}>\gamma, \ \forall y\in \mathbb{Y}.
\end{equation}
Define $f:\mathbb{Y}\rightarrow \mathbb{R}$ as $f(y):=-|||y|||,\forall y\in \mathbb{Y}$ and let 
$$
A_1:=\{0_\mathbb{Y}\}\times (-\infty, 0], A_2:=\epi f\ \ {\rm and} \ \ \bar x:=(0_{\mathbb{Y}}, 0).
$$
We first prove that 
\begin{equation}\label{3.06}
  \mathbf{d}(x, A_1\cap A_2)\leq 2(\mathbf{d}(x, A_1)+\mathbf{d}(x, A_2)),\ \ \forall x=(y,\alpha)\in \mathbb{Y}\times \mathbb{R},
\end{equation}
which means that $\{A_1, A_2\}$ is substransversal.

Indeed, for any $x=(y,\alpha)\in \mathbb{Y}\times \mathbb{R}$, one has
\begin{equation}\label{3.07}
  \mathbf{d}(x, A_1\cap A_2)=\|y\|+|\alpha|\ \ {\rm and} \ \ \mathbf{d}(x, A_1)=\left\{
  \begin{aligned}
  \|y\|, \ \ \ & \alpha<0,\\
   \|y\|+\alpha, \  & \alpha\geq 0.
  \end{aligned}
  \right.
\end{equation}

If $\alpha\geq 0$, then \eqref{3.06} holds since $\mathbf{d}(x, A_1\cap A_2)=\mathbf{d}(x, A_1)$.

If $-|||y|||\leq \alpha <0$, one has
$$
\mathbf{d}(x, A_1\cap A_2)=\|y\|+|\alpha|\leq 2\|y\|= 2\mathbf{d}(x, A_1)
$$
and thus \eqref{3.06} holds.

If $\alpha <-|||y|||$, then for any $(z,\mu)\in A_2$, one has $-|||z|||\leq \mu$ and 
\begin{eqnarray*}
-|||y|||-\alpha&\leq& -|||y|||-\alpha+|||z|||+\mu\\
&\leq&\big||||y|||-|||z|||\big|+|\alpha-\mu|\\
&\leq& |||y-z|||+|\alpha-\mu|\\
&\leq& \|y-z\|+|\alpha-\mu|,
\end{eqnarray*}
the last inequality following  by \eqref{3.05}. This implies that
$$
-|||y|||-\alpha\leq \mathbf{d}(x, A_2)
$$
and then \eqref{3.05} and \eqref{3.07} give
$$
\mathbf{d}(x, A_1\cap A_2)=\|y\|-\alpha=\|y\|+|||y|||-|||y|||-\alpha\leq 2\|y\|+\mathbf{d}(x, A_2).
$$
Hence \eqref{3.06} holds.

Next, we show that for any $\epsilon>0$, one has
\begin{equation}\label{3.07a}
 \hat N(A_1\cap A_2, \bar x)\not\subseteq\bigcup\Big\{\hat N(A_1, x_1)+\hat N(A_2, x_2): x_i\in A_i\cap \mathbf{B}(\bar x, \varepsilon),\ i=1,2 \Big\}+\varepsilon \mathbf{B}_{\mathbb{X}^\ast },
\end{equation}
which contradicts (iii).

Note that  $\hat N(A_1\cap A_2,\bar x)=\mathbb{Y}^\ast \times \mathbb{R}$. Take any $\epsilon>0$. For any $x_1=(y_1,\alpha_1)\in A_1\cap \mathbf{B}(\bar x,\epsilon)$, one has $y_1=0_\mathbb{Y}$ and
$$
\hat N(A_1,x_1)=\left\{
\begin{aligned}
  \mathbb{Y}^\ast \times\{0\},\ \ \ \ \ \ & \alpha_1<0, \\
  \mathbb{Y}^\ast \times [0,+\infty), \ & \alpha_1=0.
\end{aligned}
\right.
$$
\textbf{Claim}. Let  $x_2=(y_2,\alpha_2)\in A_2\cap \mathbf{B}(\bar x,\epsilon)$.  Then %
\begin{equation}\label{3.08}
 \hat N(A_2,x_2)=\{(0,0)\}.
\end{equation}
Let $(y^\ast , -\lambda)\in\hat N(A_2,x_2)$. Suppose on the contrary that $\lambda\not=0$. Applying \cref{lem2.2a}, we have $\lambda>0$ and $\frac{y^\ast }{\lambda}\in\hat\partial f(y_2)$. Thus,
$$
\liminf_{\|h\|\rightarrow 0}\frac{-|||y_2+h|||+|||y_2|||-\langle\frac{y^\ast }{\lambda}, h\rangle}{\|h\|}\geq 0.
$$
Hence
$$
\limsup_{\|h\|\rightarrow 0}\frac{|||y_2+h|||+|||y_2-h|||-2|||y_2|||}{\|h\|}\leq 0,
$$
which is a contradiction with \eqref{3.05}. The claim  is established.
\vskip 2mm
By the definition, for any $\epsilon>0$, there exists $\delta>0$ such that
$$
\langle y^\ast , y-y_2\rangle\leq\epsilon(\|y-y_2\|+|\alpha-\alpha_2|),\ \ \forall (y,\alpha)\in (y_2,\alpha_2)+\delta \mathbf{B}_{\mathbb{Y}\times\mathbb{R}}\cap A_2.
$$
For each $y\in y_2+\frac{\delta}{2} \mathbf{B}_\mathbb{Y}$, take $\alpha:=-|||y|||+|||y_2|||+\alpha_2$. Then one has 
$$
\|(y,\alpha)-(y_2,\alpha_2)\|=\|y-y_2\|+|\alpha-\alpha_2|\leq 2\|y-y_2\|
$$
and from \eqref{3.05} we get that 
\begin{eqnarray*}
\langle y^\ast , y-y_2\rangle\leq\epsilon(\|y-y_2\|+|\alpha-\alpha_2|)\leq 2\epsilon\|y-y_2\|.
\end{eqnarray*}
This implies that $y^\ast =0$ and thus \eqref{3.08} holds.

Noticing that  $\hat N(A_1\cap A_2,\bar x)=\mathbb{Y}^\ast \times \mathbb{R}$ and $\varepsilon\mathbf{B}_{\mathbb{X}^\ast }\subseteq \mathbb{Y}^\ast \times [-\epsilon, \epsilon]$, this yields
$$
\hat N(A_1\cap A_2, \bar x)\not\subset \mathbb{Y}^\ast \times [0,+\infty)+\{(0,0)\}+\varepsilon \mathbf{B}_{X^\ast },
$$
and therefore  it follows that \eqref{3.07a} holds. 
The proof is complete. \hfill$\Box$

\begin{remark}\label{rem3.1}

	 It would be ideal if  in (iii) of \cref{th3.1}, $\epsilon$  could  be taken as $0$ as this consequence would imply the strong Fr\'echet CHIP of $\{A_1,\cdots, A_m\}$ at $\bar x$. However, the following example shows that such ideal case may not hold even in an Euclidean space. 

\end{remark}

\begin{example} Let $\mathbb{X}:=\mathbb{R}^2$, $A_1:=\mathbb{R}\times\{0\}\cup\{(x_1,x_2)\in \mathbb{X}: x_1=x_2\}$, $A_2:=\{0\}\times (-\infty, 0]$ and $\bar x:=(0,0)$. Then $A_1\cap A_2=\{\bar x\}$ and one can verify that $\{A_1,A_2\}$ is subtransversal at $\bar x$. However, by the computation, one has
	$$
	\hat N(A_1\cap A_2, \bar x)=\mathbb{R}^2, \hat N(A_1, \bar x)=\{(0,0)\}\ \ {\rm and} \ \ \hat N(A_2, \bar x)=\mathbb{R}\times (-\infty, 0],
	$$
	which implies that
	$$
	\hat N(A_1\cap A_2, \bar x)\not\subseteq\hat N(A_1, \bar x)+\hat N(A_2, \bar x).
	$$
	Hence  $\{A_1,A_2\}$ does not have the strong Fr\'echet CHIP at $\bar x$.
\end{example}

\fcolorbox{black}{yellow}{
\begin{minipage}{12cm}
\textit{In the remainder of this paper, unless otherwise stated, we suppose that $\mathbb{X}$ is  Asplund and}
$\{A_1,\cdots,A_m\}$ \textit{is a collection of finitely many closed sets in $\mathbb{X}$ such that the intersection is nonempty}.
\end{minipage}
}
 
\medskip

We  are now in a position to consider circumstances when one can pass to the limit in \eqref{3.11} as $\varepsilon\downarrow 0$ and thus obtain necessary conditions for subtransversality of $\{A_1,\cdots,A_m\}$ in terms of limiting normal cones. The following theorem shows that it can be done in the finite-dimensional setting.

\begin{theorem}\label{th3.3}
Let $\bar x\in A:=\bigcap_{i=1}^mA_i$. Suppose that $\mathbb{X}$ is finite-dimensional and $\{A_1,\cdots,A_m\}$ is  subtransversal at $\bar x$. Then there exist $\tau,\delta\in (0, +\infty)$ such that
\begin{equation}\label{3.14}
    N(A,x)\cap \mathbf{B}_{\mathbb{X}^\ast }\subseteq\tau\sum_{i=1}^m N(A_i, x)\cap \mathbf{B}_{\mathbb{X}^\ast }
\end{equation}
holds for all $x\in \mathbf{B}(\bar x,\delta)\cap A$.
\end{theorem}

{\bf Proof.} Since $\{A_1,\cdots,A_m\}$ is  subtransversal at $\bar x$, then there exist $\tau, r\in(0,+\infty)$ such that \eqref{3.3} holds. Let $\delta:=\frac{r}{2}$ and $x\in A\cap \mathbf{B}(\bar x,\delta)$. Take any $x^\ast \in N(A,x)\cap \mathbf{B}_{\mathbb{X}^\ast }$. Since $\mathbb{X}$ is finite-dimensional, it follows that there exist $x_k\stackrel{A}\rightarrow x$ and $x^\ast _k\rightarrow x^\ast $ such that $x^\ast _k\in\hat N(A,x_k)$ for all $k$. For any $k$, set $\widetilde{x}^\ast _k:=\frac{x_k^\ast }{\|x_k^\ast \|}$. Similarly to the proof of  (i)\;$\Rightarrow$\;(ii) in \cref{th3.1}, one can get
$$
\widetilde{x}_k^\ast \in\hat\partial \Big(\tau\sum_{i=1}^m\mathbf{d}(\cdot,A_i)+\frac{1}{k}\|\cdot-x_k\|\Big)(x_k).
$$
By using \cref{lem2.1,lem2.2}, there are $u_{i,k}\in \mathbf{B}(x_k,\frac{1}{k})$, $u_{i,k}^\ast \in\hat\partial \mathbf{d}(\cdot, A_i)(u_{i,k})$, $x_{i,k}\in A_i$ and $x_{i,k}^\ast \in\hat N(A_i, x_{i,k}), i=1,\cdots,m,$ such that
\begin{equation}\label{3.15}
  \|x_{i,k}-u_{i,k}\|<\mathbf{d}(u_{i,k}, A_i)+\frac{1}{k},\ \ \|x^\ast _{i,k}-u^\ast _{i,k}\|<\frac{1}{k}
\end{equation}
and
\begin{equation}\label{3.16}
  \widetilde{x}^\ast _k\in\tau\sum_{i=1}^mu_{i,k}^\ast +\frac{1}{k}\mathbf{B}_{\mathbb{X}^\ast }+\frac{1}{k}\mathbf{B}_{\mathbb{X}^\ast }\subset\tau\sum_{i=1}^mx_{i,k}^\ast +\frac{\tau m+2}{k}\mathbf{B}_{\mathbb{X}^\ast }.
\end{equation}
For any $i$, since $\{\|x_{i,k}^\ast \|\}_{k=1}^{\infty}$ is bounded and $\mathbb{X}$ is finite-dimensional, we can assume that $x_{i,k}^\ast \rightarrow x_i^\ast $ as $k\rightarrow \infty$ (considering subsequence if necessary) and it is not hard to verify that $x_{i,k}\stackrel{A_i}\rightarrow x$ as $k\rightarrow \infty$ (thanks to \eqref{3.15} and $x_k\rightarrow x$). From this, it follows that $x_i^\ast \in N(A_i, x)$ and
$$
\|x_i^\ast \|=\lim_k\|x^\ast _{i,k}\|\leq\limsup_k\|u^\ast _{i,k}\|\leq 1.
$$
Taking limits as $k\rightarrow \infty$ in \eqref{3.16}, one has
$$
\frac{x^\ast }{\|x^\ast \|}=\tau\sum_{i=1}^m
x^\ast _i.
$$
This implies that
$$
x^\ast =\tau\sum_{i=1}^m
\|x^\ast \|x^\ast _i\subset \tau\sum_{i=1}^m
N(A_i,x)\cap \mathbf{B}_{\mathbb{X}^\ast }
$$
as $\|x^\ast \|\leq 1$ and $\|x^\ast _i\|\leq 1$. Hence \eqref{3.14} holds. The proof is complete. \hfill$\Box$

\medskip


\begin{remark}\label{re3.3}
	(a) It is noted that \eqref{3.14} may not imply property (G) of limiting normal cones $\{N(A_1,x),\cdots,N(A_m,x)\}$ necessarily since $N(A,x)\supseteq\sum_{i=1}^mN(A_i,x)$ may not hold trivially. Further, for the case that $A_1,\cdots, A_m$ are closed convex sets, it is easy to verify that \eqref{3.14} holds if and only if  $\{A_1,\cdots,A_m\}$ has the strong CHIP at $x$ and $\{N(A_1,x),\cdots,N(A_m,x)\}$ has property (G) with constant $\tau>0$. This and {\bf Theorem A} imply that subtransversality is equivalent to \eqref{3.14} which holds for all $x\in A$ close to $\bar x$. 

(b) It is known from \cref{th3.3} that the subtransversality of $\{A_1,\cdots,A_m\}$ implies the validity of \eqref{3.14} around $\bar x$. However,  the inverse implication may not hold necessarily. Consider {\it Example 2.1}(b). 
Then $A_1\cap A_2=\{\bar x\}$ and by computation, one has $N(A_1\cap A_2,\bar x)=\mathbb{X}$,
$$
N(A_2,\bar x)=\mathbb{R}_-\times\mathbb{R}_- \ \ {\rm and} \ \ N(A_1,\bar x) = ([0, +\infty)\times\{0\})\cup(\{0\}\times (-\infty, 0])\cup\{(t,t):t\in\mathbb{R}\}.
$$
Thus, one can verify that there exists $\tau>0$ such that \eqref{3.14} holds with $x=\bar x$.

On the other hand, for each $k\in \mathbb{N}$, let $x_k=\big(\frac{1}{k}, (\frac{2}{k}-\frac{1}{k^2})^{\frac{1}{2}}\big)$. By computation, one has
$$
\mathbf{d}(x_k, A_1\cap A_2)=\big(\frac{2}{k}\big)^{\frac{1}{2}}, \mathbf{d}(x_k, A_1)=\frac{1}{k}\,\,\,{\rm and}\,\,\,\mathbf{d}(x_k, A_2)=0.
$$
Thus
$$
\frac{\mathbf{d}(x_k, A_1\cap A_2)}{\mathbf{d}(x_k, A_1)+\mathbf{d}(x_k, A_2)}\rightarrow \infty \ ({\rm as}\ k\rightarrow \infty),
$$
which means that $\{A_1,A_2\}$ is not  subtransversal at $\bar x$.

\end{remark}




Further, we consider the subtransversality for a collection of some special closed sets in convex-composite optimization. It is proved that the equivalence result on subtransversality, strong Fr\'echet CHIP and property (G) still holds for this case. Before doing this, we need the following lemma which is of independent interest.


\begin{lemma}\label{lem3.1}
	Let $f:X\rightarrow Y$ be a continuously differentiable mapping between Banach spaces,  $C$ be a closed convex set in $Y$ and $\bar x\in A:=f^{-1}(C)$ such that $\triangledown f(\bar x)$ is surjective. Then there exist $\ell,L,r\in (0,+\infty)$ such that
	\begin{equation}\label{3.017}
		\triangledown f(x)^\ast (Y^\ast )\cap\ell \mathbf{B}_{X^\ast }\subseteq \triangledown f(x)^\ast (\mathbf{B}_{Y^\ast })\subseteq L\mathbf{B}_{X^\ast },\ \ \forall x\in \mathbf{B}(\bar x,r)
	\end{equation}
	and
	\begin{equation}\label{3.017a}
		\hat N(A,x)\cap \ell\mathbf{B}_{X^\ast }\subseteq \triangledown f(x)^\ast (N(C,f(x))\cap \mathbf{B}_{Y^\ast })\subseteq \hat N(A,x)\cap L\mathbf{B}_{X^\ast }, \ \ \forall x\in \mathbf{B}(\bar x,r)\cap A,
	\end{equation}
where $\triangledown f(x)^*$ is the adjoint of  $\triangledown f(x)$.
\end{lemma}

{\bf Proof.} We first show that there exist $\ell>0$ and $r\in (0, r_0)$ such that
\begin{equation}\label{3-18}
	\|\triangledown f(x)-\triangledown f(\bar x)\|<\ell\ \ {\rm and} \ \ \ell \mathbf{B}_{Y}\subseteq\triangledown f(x)(\mathbf{B}_{X})\ \ \forall x\in \mathbf{B}(\bar x, r).
\end{equation}

Applying the celebrated  Open Mapping Theorem, $0_Y$ is an interior point of $\triangledown f(\bar x)(\mathbf{B}_{X})$ and thus there exists $\ell >0$ such that
\begin{equation}\label{3-19}
	3\ell \mathbf{B}_{Y}\subseteq\triangledown f(\bar x)(\mathbf{B}_{X}).
\end{equation}
Noting that $f$ is continuously differentiable at $\bar x$, then there exists $r>0$ such that $f$ is Fr\'{e}chet differentiable on $\mathbf{B}(\bar x,r)$ and
\begin{equation}\label{3-20}
	\|\triangledown f(x)-\triangledown f(\bar x)\|<\ell \ \ \forall x\in \mathbf{B}(\bar x,r).
\end{equation}
Then by \eqref{3-19} and \eqref{3-20}, one has
$$
3\ell \mathbf{B}_{Y}\subseteq(\triangledown f(x)+\triangledown f(\bar x)-\triangledown f(x))(\mathbf{B}_{X})\subseteq \triangledown f(x)(\mathbf{B}_{X})+\ell \mathbf{B}_{Y}\ \ \forall x\in \mathbf{B}(\bar x,r).
$$
By the R{\aa}dstr\"om Cancellation Lemma (cf. \cite[Lemma 2.3]{Ra}), one has
\begin{equation}
	2 \ell \mathbf{B}_{Y}\subseteq{\rm cl}(\triangledown f(x)(\mathbf{B}_{X}))\ \ \forall x\in \mathbf{B}(\bar x,r).
\end{equation}
Note that $\triangledown f(x)(\mathbf{B}_{X})$ is convex and then \cite[P.183, Theorem A.1]{Jameson} implies that $\triangledown f(x)(\mathbf{B}_{X})$ and ${\rm cl}(\triangledown f(x)(\mathbf{B}_{X}))$ have the same interior and consequently
$$
\ell \mathbf{B}_{Y}\subseteq\triangledown f(x)(\mathbf{B}_{X})\ \ \forall x\in \mathbf{B}(\bar x, r).
$$
Hence \eqref{3-18} holds.

By the continuity of $\triangledown f$, without loss of generality we can assume that there is an $L>0$ such that
\begin{equation*}
	\|\triangledown f(x)-\triangledown f(\bar x)\|<\ell\ \ {\rm and} \ \ \|\triangledown f(x)\|\leq L,\ \ \forall x\in \mathbf{B}(\bar x, r)
\end{equation*}
(taking a smaller $r$ if necessary). 

We first prove \eqref{3.017}. Let $x\in \mathbf{B}(\bar x,r)$ and $x^\ast \in\triangledown f(x)^\ast (\mathbb{Y}^\ast )\cap\ell \mathbf{B}_{X^\ast }$. Then there exists $y^\ast \in \mathbb{Y}^\ast $ such that $x^\ast =\triangledown f(x)^\ast (y^\ast )$. By virtue of \eqref{3-18}, one has
\begin{eqnarray*}
	\|y^\ast \|=\sup_{y\in \mathbf{B}_{Y}}\langle y^\ast ,y\rangle\leq\sup_{y\in\frac{1}{\ell}\triangledown f(x)(\mathbf{B}_{X})}\langle y^\ast ,y\rangle&=&\frac{1}{\ell}\sup_{u\in \mathbf{B}_{X}}\langle y^\ast ,\triangledown f(x)(u)\rangle\\
	&=&\frac{1}{\ell }\sup_{u\in \mathbf{B}_{X}}\langle x^\ast ,u\rangle=\frac{\|x^\ast \|}{\ell}\leq 1.
\end{eqnarray*}
Hence $x^\ast \in \triangledown f(x)^\ast (\mathbf{B}_{\mathbb{Y}^\ast })$. Note that $\|\triangledown f(x)\|\leq L$ and thus $\triangledown f(x)^\ast (\mathbf{B}_{\mathbb{Y}^\ast })\subseteq L\mathbf{B}_{\mathbb{X}^\ast }$. This means that \eqref{3.017} holds.

We next prove \eqref{3.017a}. Let $x\in \mathbf{B}(\bar x,r)\cap A$ and $x^\ast \in \hat N(A_i,x)\cap \ell\mathbf{B}_{\mathbb{X}^\ast }$. Note that $A=f^{-1}(C)$ and it follows from \cite[Corollary 1.15]{Mordukhovich} that
\begin{equation}\label{3-26}
	\hat N(A, x)=\triangledown f(x)^\ast (N(C,f(x))).
\end{equation}
Then there exists $y_1^\ast \in N(C,f(x))$ such that $x^\ast =\triangledown f(x)^\ast (y_1^\ast )$. This and \eqref{3-23} imply that there exists $y_2^\ast \in \mathbf{B}_{\mathbb{Y}^\ast }$ such that $x^\ast =\triangledown f(x)^\ast (y_2^\ast )$. Since $\triangledown f(x)^\ast $ is one-to-one (due to \eqref{3-18}), it follows that $y_1^\ast =y_2^\ast $ and thus $x^\ast \in\triangledown f(x)^\ast (N(C,f(x))\cap \mathbf{B}_{\mathbb{Y}^\ast })$. Note that $\|\triangledown f(x)\|\leq L$ and so \eqref{3.017a} follows immediately from \eqref{3-26}. The proof is complete. \hfill$\Box$

\medskip

The following theorem shows that the subtransversality is characterized by  the strong Fr\'echet CHIP and property (G) for some special closed sets in convex-composite optimization in the Asplund space setting.

\begin{theorem}\label{th3.4}
Let $\mathbb{Y}$ be an Asplund space, $f:\mathbb{X}\rightarrow \mathbb{Y}$ be a continuously differentiable mapping and  $\{C_1,\cdots,C_m\}$ be a collection of finitely many closed convex sets in $\mathbb{Y}$ such that $C:=\bigcap_{i=1}^mC_i$ is nonempty. Let $A_i:=f^{-1}(C_i)(i=1,\cdots,m)$ and $\bar x\in A:=\bigcap_{i=1}^mA_i$ such that $\triangledown f(\bar x)$ is surjective. Then the following statements are equivalent:
\begin{itemize}
\item[\rm (i)] $\{A_1,\cdots,A_m\}$ is  subtransversal at $\bar x$;
\item[\rm (ii)] there exist $\tau,\delta>0$ such that for all $x\in A\cap \mathbf{B}(\bar x,\delta)$, one has that $\{A_1,\cdots,A_m\}$ has the strong Fr\'echet CHIP and $\{\hat N(A_1,x),\cdots,\hat N(A_m,x)\}$ has property (G) with constant $\tau$;
\item[\rm (iii)] there exist $\tau,\delta>0$ such that for all $x\in A\cap \mathbf{B}(\bar x,\delta)$, one has that $\{A_1,\cdots,A_m\}$ has the strong limiting CHIP and $\{N(A_1,x),\cdots,N(A_m,x)\}$ has property (G) with constant $\tau$.
\end{itemize}
\end{theorem}

{\bf Proof.} Since $\triangledown f(\bar x)$ is surjective, then \cite[Theorem 1.57]{Mordukhovich} implies that $f$ is metrically regular at $\bar x$ and thus there exist $\kappa,r_0>0$ such that
\begin{equation}\label{3-17}
  \mathbf{d}(x, f^{-1}(y))\leq\kappa \|y-f(x)\|\ \ \forall (x,y)\in \mathbf{B}(\bar x,r_0)\times \mathbf{B}(f(\bar x), r_0).
\end{equation}
Applying \cref{lem3.1}, there exist $\ell,L>0$ and $r_1 \in (0,r_0)$ such that
\begin{equation}\label{3-23}
  \triangledown f(x)^\ast (\mathbb{Y}^\ast )\cap\ell \mathbf{B}_{\mathbb{X}^\ast }\subseteq \triangledown f(x)^\ast (\mathbf{B}_{\mathbb{Y}^\ast })\subseteq L\mathbf{B}_{\mathbb{X}^\ast },\ \ \forall x\in \mathbf{B}(\bar x,r_1),
\end{equation}
\begin{equation}\label{3-24}
  \hat N(A,x)\cap \ell\mathbf{B}_{\mathbb{X}^\ast }\subseteq \triangledown f(x)^\ast (N(C,f(x))\cap \mathbf{B}_{\mathbb{Y}^\ast })\subseteq \hat N(A,x)\cap \ell\mathbf{B}_{\mathbb{X}^\ast }
\end{equation}
and
\begin{equation}\label{3-25}
  \hat N(A_i,x)\cap \ell\mathbf{B}_{\mathbb{X}^\ast }\subseteq \triangledown f(x)^\ast (N(C_i,f(x))\cap \mathbf{B}_{\mathbb{Y}^\ast })\subseteq \hat N(A_i,x)\cap L\mathbf{B}_{\mathbb{X}^\ast }
\end{equation}
hold for all $x\in \mathbf{B}(\bar x,r_1)\cap A$ and for each $i=1,\cdots,m$.  

Now, we first prove the equivalence between (i) and (ii). By virtue of \cite[Corollary 4.2]{ZW}, one has that $\{A_1,\cdots,A_m\}$ is  subtransversal at $\bar x$ if and only if $\{C_1,\cdots,C_m\}$ is  subtransversal at $f(\bar x)$. Based on {\bf Theorem A} in  the Introduction, it suffices to prove that (ii) holds if and only if there exist $\tau_1,\delta_1>0$ such that
\begin{equation}\label{3-27}
  N(C,y)\cap \mathbf{B}_{\mathbb{Y}^\ast }\subset \tau_1\sum_{i=1}^m(N(C_i,y)\cap \mathbf{B}_{\mathbb{Y}^\ast })
\end{equation}
holds for all $y\in C\cap \mathbf{B}(f(\bar x),\delta_1)$.

Suppose that (ii) holds. Then one can verify that
\begin{equation}\label{3-28}
  \hat N(A,x)\cap \mathbf{B}_{\mathbb{X}^\ast }\subset\tau\sum_{i=1}^m(\hat N(A_i,x)\cap \mathbf{B}_{\mathbb{X}^\ast }),\ \ \forall x\in \mathbf{B}(\bar x,\delta).
\end{equation}
Without loss  of generality, we can assume that $\delta\in(0,r_1)$ (taking a smaller $\delta$ if necessary). Take $\delta_1\in (0,\delta)$ such that $\kappa\delta_1<\delta$. Let $y\in C\cap \mathbf{B}(f(\bar x),\delta_1)$. Then by \eqref{3-17}, one has
$$
\mathbf{d}(\bar x, f^{-1}(y))\leq\kappa\|f(\bar x)-y\|<\kappa\delta_1
$$
and thus there is $x\in f^{-1}(y)$ such that $\|x-\bar x\|<\kappa\delta_1$. This implies that $x\in A\cap \mathbf{B}(\bar x,\delta)$. By virtue of \eqref{3-24}, \eqref{3-25} and \eqref{3-28}, one has
\begin{eqnarray*}
\triangledown f(x)^\ast (N(C,y)\cap \mathbf{B}_{\mathbb{Y}^\ast })&=&\triangledown f(x)^\ast (N(C,f(x))\cap \mathbf{B}_{\mathbb{Y}^\ast })\subset\hat N(A,x)\cap L\mathbf{B}_{\mathbb{X}^\ast }\\&\subset&L\tau \sum_{i=1}^m(\hat N(A_i,x)\cap \mathbf{B}_{\mathbb{X}^\ast })\\
&=&\frac{L\tau}{\ell} \sum_{i=1}^m(\hat N(A_i,x)\cap \ell\mathbf{B}_{\mathbb{X}^\ast })\\
&\subset&\frac{L\tau}{\ell} \sum_{i=1}^m \triangledown f(x)^\ast (N(C_i,f(x))\cap \mathbf{B}_{\mathbb{Y}^\ast })\\
&=&\triangledown f(x)^\ast \left(\frac{L\tau}{\ell} \sum_{i=1}^m(N(C_i,f(x))\cap \mathbf{B}_{\mathbb{Y}^\ast })\right).
\end{eqnarray*}
Since $\triangledown f(x)^\ast $ is one-to-one, it follows that
$$
N(C,y)\cap \mathbf{B}_{\mathbb{Y}^\ast }\subset \frac{L\tau}{\ell} \sum_{i=1}^m(N(C_i, y)\cap \mathbf{B}_{\mathbb{Y}^\ast }).
$$
This means that \eqref{3-27} holds with $\tau_1:=\frac{L\tau}{l}$.

Conversely, suppose that there exist $\tau_1,\delta_1>0$ such that \eqref{3-27} holds. By the continuity of $f$, there exists $\delta\in(0, r_1)$ such that
\begin{equation}\label{3-29}
  \|f(x)-f(\bar x)\|<\delta_1,\ \ \forall x\in \mathbf{B}(\bar x,\delta).
\end{equation}
Then, for any $x\in \mathbf{B}(\bar x,\delta)\cap A$, by \eqref{3-29}, one has $f(x)\in \mathbf{B}(f(\bar x),\delta_1)$ and it follows from \eqref{3-24}, \eqref{3-25} and \eqref{3-26} that
\begin{eqnarray*}
\hat N(A,x)\cap \mathbf{B}_{\mathbb{X}^\ast }&\subset&\frac{1}{\ell}\triangledown f(x)^\ast (N(C,f(x))\cap \mathbf{B}_{\mathbb{Y}^\ast })\subset\frac{1}{\ell}\triangledown f(x)^\ast \left(\tau_1 \sum_{i=1}^mN(C_i,f(x))\cap \mathbf{B}_{\mathbb{Y}^\ast }\right)\\&=&\frac{\tau_1}{\ell}\sum_{i=1}^m\triangledown f(x)^\ast (N(C_i,f(x))\cap \mathbf{B}_{\mathbb{Y}^\ast })\\
&\subset&\frac{\tau_1}{\ell}\sum_{i=1}^m(\hat N(A_i,x)\cap L\mathbf{B}_{\mathbb{X}^\ast })\\
&=&\frac{L\tau_1}{\ell}\sum_{i=1}^m(\hat N(A_i,x)\cap \mathbf{B}_{\mathbb{X}^\ast }).
\end{eqnarray*}
Hence (ii) holds with $\tau:=\frac{L\tau_1}{\ell}$ and $\delta>0$.

We next prove the equivalence between (ii) and (iii). Suppose that (ii) holds. Without loss of generality, we can assume that $\delta\in (0,r_1)$ (taking a smaller $\delta$ if necessary). We claim that
\begin{equation}\label{3-30}
  N(A,x)=\hat N(A, x)\ \ {\rm and}\ \ N(A_i,x)=\hat N(A_i, x)\ \  \forall x\in \mathbf{B}(\bar x,\delta)\cap A \ \ {\rm and} \ \ \forall i=1,\cdots,m.
\end{equation}

Let $x\in \mathbf{B}(\bar x,\delta)\cap A$ and $x^\ast \in N(A,x)$. Then there exist $x_k\stackrel{A}\rightarrow x$ and $x_k^\ast \stackrel{w^\ast }\rightarrow x^\ast $ such that $x_k^\ast \in \hat N(A,x_k)$ for all $k$. Applying the  Banach-Steinhaus Theorem, 
 we can assume that $\|x_k^\ast \|\leq M$ for each $k$ and some $M>0$. By \cite[Corollary 1.15]{Mordukhovich}, for any $k$ sufficiently large, there exists $y_k^\ast \in N(C, f(x_k))$ such that
$$
x^\ast _k=\triangledown f(x)^\ast (y_k^\ast ).
$$
This and \eqref{3-23} imply that
\begin{eqnarray*}
\frac{\ell}{M}x_k^\ast \in \triangledown f(x)^\ast (\mathbb{Y}^\ast )\cap \ell\mathbf{B}_{\mathbb{X}^\ast }\subset  \triangledown f(x)^\ast (\mathbf{B}_{\mathbb{Y}^\ast })
\end{eqnarray*}
and thus there is $z_k^\ast \in \mathbf{B}_{\mathbb{Y}^\ast }$ such that
$$
\frac{\ell}{M}x^\ast _k=\triangledown f(x)^\ast (z_k^\ast ).
$$
Note that $\triangledown f(x)^\ast $ is one-to-one and thus
$$
y_k^\ast =\frac{M}{\ell}z_k^\ast ,
$$
which implies that $\{\|y_k^\ast \|\}$ is bounded. Since $Y$ is an Asplund space, it follows that $\mathbf{B}_{\mathbb{Y}^\ast }$ is weak$^\ast $-sequentially compact and without loss of generality, we can assume that
$$
y_k^\ast \stackrel{w^\ast }\rightarrow y^\ast \in \mathbb{Y}^\ast 
$$
(taking subsequence if necessary). Then by the continuity of $f$ and $\triangledown f$, it is not hard to verify that $y^\ast \in N(C, f(x))$ and
$$
\triangledown f(x_k)^\ast (y_k^\ast )\stackrel{w^\ast }\rightarrow \triangledown f(x)^\ast (y^\ast )
$$
(thanks to $x_k\stackrel{A}\rightarrow x$). Noting that  $\triangledown f(x_k)^\ast (y_k^\ast )=x^\ast _k\stackrel{w^\ast }\rightarrow x^\ast $, it follows that $x^\ast =\triangledown f(x)^\ast (y^\ast )$ and consequently
$$
x^\ast \in \triangledown f(x)^\ast (N(C, f(x)))=\hat N(A,x)
$$
as $\triangledown f(x)^\ast (N(C, f(x)))=\hat N(A,x)$ follows from \cite[Corollary 1.15]{Mordukhovich}. This means that $N(A,x)=\hat N(A,x)$ and then \eqref{3-30} holds as $N(A_i,x)=\hat N(A_i,x)$ follows from the similar proof. The proof is complete. \hfill$\Box$

\begin{remark}
	It should be noted that properties (ii) and (iii) in \cref{th3.4} are essentially the same due to \eqref{3-30}. Further, for the case where  $\{A_1,\cdots,A_m\}$ is a collection of closed convex sets, \cref{th3.4} reduces to Theorem A in the Introduction.  Hence our work actually extends equivalence results among subtransversality, strong CHIP and property (G) from the convex case to the convex-composite one.
\end{remark}

\setcounter{equation}{0}
\section{Applications to Error Bounds of   Inequality Systems}

In this section, we apply results from  Section 3 to error bounds of inequality  systems. Let $\varphi_1,\cdots,\varphi_m$ be continuous functions on $\mathbb{X}$. For each $i$, let 
$$
A_i:=\{x\in X: \varphi_i(x)\leq 0\}.
$$ 
We consider the following  inequality  system:
\begin{equation}\label{4.1}
  \varphi_1(x)\leq 0, \cdots, \varphi_m(x)\leq 0.
\end{equation}
Let $A_i:=\{x\in X: \varphi_i(x)\leq 0\} $ and  $A=\bigcap_{i=1}^mA_i$ the solution set of the inequality  system \eqref{4.1}. We assume that such system has at least one solution.  For any $x\in A$, set
\begin{equation}\label{4.2a}
  I(x):=\{i:\varphi_i(x)=0\}.
\end{equation}

We first recall the definition of error bounds of inequality  systems.  

\begin{definition}\label{def4.1}
  Let $\bar x\in A$. We say that 
  \begin{itemize}
    \item [\rm (i)] the   inequality system \eqref{4.1} has a local error bound at $\bar x$, if there exists $\eta,\delta\in (0,+\infty)$ such that 
\begin{equation}\label{4.2}
  \mathbf{d}(x, A)\leq \eta\sum_{i=1}^{m}\max\{\varphi_i(x),0\},\ \ \forall x\in \mathbf{B}(\bar x,\delta);
\end{equation}
 \item [\rm (ii)] the inequality system \eqref{4.1} has a global error bound, if there exists $\eta\in (0,+\infty)$ such that 
\begin{equation*}\label{4.4a}
  \mathbf{d}(x, A)\leq \eta\sum_{i=1}^{m}\max\{\varphi_i(x),0\},\ \ \forall x\in \mathbb{X}.
\end{equation*}
  \end{itemize}

\end{definition}

\begin{remark}\label{rem4.1}
 It is noted that there is no necessary relation between the local error bound of  the   inequality system \eqref{4.1} and the subtransversality of $\{A_1,\cdots,A_m\}$. For example, if  the inequality system \eqref{4.1} has only one inequality, say $\varphi_1(x)\leq 0$, then the subtransversality (of $\{A_1\}$) is valid automatically while the local error bound of $\varphi_1$ may not hold necessarily.

\end{remark}

Next, we consider error bounds for  convex-composite inequality  systems. 

\medskip

{\it Throughout the remainder of this section, we suppose that $\mathbb{Y}$ is an Asplund space, $f:\mathbb{X}\rightarrow \mathbb{Y}$ is continuously differentiable and each $g_i:\mathbb{Y}\rightarrow \mathbb{R}$ is a continuous convex function on $\mathbb{Y}$ for $i=1,\cdots,m$}. 

\medskip

For each $i$, let $\varphi_i:=g_i\circ f$, $A_i:=\{x\in \mathbb{X}:\varphi_i(x)\leq 0\}$ and $C_i:=\{y\in \mathbb{Y}:g_i(y)\leq 0\}$. We consider the following convex-composite   inequality system:
\begin{equation}\label{4.7}
  (g_1\circ f)(x)\leq 0,\cdots, (g_m\circ f)(x)\leq 0.
\end{equation}
We still denote by $A:=\bigcap_{i=1}^mA_i$ the solution set.  It is easy to verify that each $C_i$ is closed and convex,  
\begin{equation}\label{4.8}
  A_i=f^{-1}(C_i)\ \forall i\ \ {\rm and} \ \ A=f^{-1}\left(\bigcap_{i=1}^{m} C_i\right).
\end{equation}

The following theorem provides sufficient conditions for the local error bound of convex-composite inequality  system \eqref{4.7}.
\begin{theorem}\label{th4.1}
Let $\bar x\in A$  be such that $\triangledown f(\bar x)$ is surjective, each convex inequality $g_i(y)\leq 0$ has a local error bound at $f(\bar x)$ and $\{A_1,\cdots, A_m\}$ is subtransversal  at $\bar x$. Then the  convex-composite inequality  system \eqref{4.7} has a local error bound at $\bar x$.
\end{theorem}

{\bf Proof} Since each convex inequality $g_i(y)\leq 0$ has a local error bound at $f(\bar x)$, there exist $\eta_1,r_1>0$ such that
\begin{equation}\label{4.9}
  \mathbf{d}(y,C_i)\leq\eta_1\max\{g_i(y),0\},\ \ \forall y\in \mathbf{B}(f(\bar x),r_1).
\end{equation}
Since $\{A_1,\cdots, A_m\}$ is subtransversal   at $\bar x$ and $f$ is continuous, there exist $\tau_1,\delta_1>0$ such that
\begin{equation}\label{4.10}
  \|f(x)-f(\bar x)\|<r_1\ \ {\rm and} \ \ \mathbf{d}(x,A)\leq\tau_1\sum_{i=1}^{m}\mathbf{d}(x,A_i)\ \ \forall x\in \mathbf{B}(\bar x,\delta_1).
\end{equation}
Since $A_i= f^{-1}(C_i)$ for each $i$ and $\{C_1,\cdots, C_m\}$ is a collection of closed convex sets, \cite[Proposition 3.5]{ZW} implies that there exist $\tau_2>0$ and $\delta_2\in (0, \delta_1)$ such that
\begin{equation}\label{4.11}
  \mathbf{d}(x, A_i)\leq\tau_2 \mathbf{d}(f(x), C_i)\ \ \forall x\in \mathbf{B}(\bar x,\delta_2) \ {\rm and} \ \forall i.
\end{equation}
By virtue of \eqref{4.9}, \eqref{4.10} and \eqref{4.11}, for any $x\in \mathbf{B}(\bar x,\delta_2)$, one has
\begin{eqnarray*}
  \mathbf{d}(x,A)\leq\tau_1\sum_{i=1}^{m}\mathbf{d}(x,A_i)&\leq& \tau_1\tau_2\sum_{i=1}^{m} \mathbf{d}(f(x), C_i)\\
  &\leq& \tau_1\tau_2\sum_{i=1}^{m}\max\{g_i(f(x)),0\}\\
  &=&\tau_1\tau_2\sum_{i=1}^{m}\max\{\varphi_i(x),0\}.
\end{eqnarray*}
This implies that the  convex-composite inequality  system \eqref{4.7} has a local error bound at $\bar x$ with constant $\eta:=\tau_1\tau_2$. The proof is complete.\hfill$\Box$

\medskip

Now, we are in a  position  to prove dual characterizations for error bounds of  the convex-composite inequality system \eqref{4.7}. Before doing this, we need the following lemma on subdifferential calculus for the convex-composite function.

\begin{lemma}\label{lem4.1}
  Let $\bar x\in \mathbb{X}$. Then 
  \begin{equation}\label{4.12}
    \hat\partial\varphi_i(\bar x)=\partial\varphi_i(\bar x)=\triangledown f(\bar x)^\ast (\partial g_i(f(\bar x))),\ \ \forall i=1,\cdots,m.
  \end{equation}
\end{lemma}

{\bf Proof.} Let $i$ be given. We first show that
\begin{equation}\label{4.13a}
   \hat\partial\varphi_i(\bar x)=\triangledown f(\bar x)^\ast (\partial g_i(f(\bar x))).
\end{equation}
Note that $g_i$ is a continuous convex function and thus $g_i$ is locally Lipschtizian around $f(\bar x)$. Then there exist $\kappa,r>0$ such that 
\begin{equation}\label{4.13}
  |g_i(y_1)-g_i(y_2)|\leq\kappa\|y_1-y_2\|, \ \ \forall y_1,y_2\in \mathbf{B}(f(\bar x),r).
\end{equation}
Let $\epsilon>0$. Since $f$ is differentiable at $\bar x$, there exists $\delta>0$ such that
\begin{equation*}
  \|f(x)-f(\bar x)\|\leq r\ \ {\rm and}\ \ \|f(x)-f(\bar x)-\triangledown f(\bar x)(x-\bar x)\|<\epsilon\|x-\bar x\|,\ \ \forall x\in \mathbf{B}(\bar x,\delta).
\end{equation*}
This and \eqref{4.13} imply that for any $x$ close to $\bar x$, one has
\begin{equation}\label{4.14}
\frac{|g_i(f(x))-g_i(f(\bar x)+\triangledown f(\bar x)(x-\bar x))|}{\|x-\bar x\|}\leq\kappa\epsilon, \ \ {\rm for \ all} \ x\ {\rm close \ to} \ \bar x.
\end{equation}

Let $y^\ast \in\partial g_i(f(\bar x))$. Then \eqref{4.14} gives that
\begin{eqnarray*}
 && \liminf_{x\rightarrow \bar x}\frac{\varphi_i(x)-\varphi_i(\bar x)-\langle \triangledown f(\bar x)^\ast (y^\ast ), x-\bar x\rangle}{\|x-\bar x\|} \\
  &=& \liminf_{x\rightarrow \bar x}\frac{g_i(f(x))-g_i(f(\bar x))-\langle y^\ast ,\triangledown f(\bar x)(x-\bar x))\rangle}{\|x-\bar x\|} \\
  &\geq&\liminf_{x\rightarrow \bar x}\frac{g_i(f(x))-g_i(f(\bar x))-g_i(f(\bar x)+\triangledown f(\bar x)(x-\bar x))+g_i(f(\bar x))}{\|x-\bar x\|} \\
  &\geq&\liminf_{x\rightarrow \bar x}\frac{g_i(f(x))-g_i(f(\bar x)+\triangledown f(\bar x)(x-\bar x))}{\|x-\bar x\|} \\
  &\geq&-\kappa\epsilon.
\end{eqnarray*}
Taking limits as $\epsilon\downarrow 0$, one has $$
\triangledown f(\bar x)^\ast (y^\ast )\in  \hat\partial\varphi_i(\bar x).
$$

On the other hand, suppose on the contrary that there exists $x^\ast \in \hat\partial\varphi_i(\bar x)\backslash \triangledown f(\bar x)^\ast (\partial g_i(f(\bar x)))$. Since $\triangledown f(\bar x)^\ast $ is weak$^\ast $-to-weak$^\ast $ continuous and $\partial g_i(f(\bar x))$ is weak$^\ast $-compact, it follows that $\triangledown f(\bar x)^\ast (\partial g_i(f(\bar x)))$ is weak$^\ast $-compact and convex. By the Separation Theorem (cf. \cite[Theorem 3.4]{Rudin}), there exist $h\in \mathbb{X}$ with  $\|h\|=1$  and $\gamma\in\mathbb{R}$ such that
\begin{equation}\label{4.15}
 \langle x^\ast , h\rangle>\gamma>\sup\{\langle x^\ast , h\rangle: x^\ast \in\triangledown f(\bar x)^\ast (\partial g_i(f(\bar x)))\}=d^+g_i(f(\bar x), \triangledown f(\bar x)h).
\end{equation}
Note that $x^\ast \in \hat\partial\varphi_i(\bar x)$ and so there exists $\delta>0$ such that
\begin{equation*}
  g_i(f(x))-g_i(f(\bar x))-\langle x^\ast , x-\bar x\rangle\geq -\epsilon\|x-\bar x\|\ \ \forall x\in \mathbf{B}(\bar x,\delta).
\end{equation*}
Then for any $t>0$ sufficiently small, one has
\begin{eqnarray*}
 \gamma-\epsilon< \langle x^\ast ,h\rangle-\epsilon&\leq& \frac{g_i(f(\bar x+th))-g_i(f(\bar x))}{t} \\&=&\frac{g_i(f(\bar x+th))-g_i(f(\bar x)+\triangledown f(\bar x)h)}{t} +\frac{g_i(f(\bar x)+\triangledown f(\bar x)h)-g_i(f(\bar x))}{t}. 
\end{eqnarray*}
Letting $t\rightarrow 0^+$ and then $\epsilon\downarrow 0$, it follows from \eqref{4.14} that
$$
d^+g_i(f(\bar x), \triangledown f(\bar x)h)\geq \gamma,
$$
which contradicts \eqref{4.15}. Hence \eqref{4.13a} holds.

We next prove that
\begin{equation}\label{4.17}
  \partial\varphi_i(\bar x)=\triangledown f(\bar x)^\ast (\partial g_i(f(\bar x))).
\end{equation}
Granting this, it follows from \eqref{4.13a} that \eqref{4.12} holds.

Note that $\partial\varphi_i(\bar x)\supseteq\triangledown f(\bar x)^\ast (\partial g_i(f(\bar x)))$ by \eqref{4.13a} and so it suffices to prove the inclusion. To this aim, let $x^\ast \in \partial\varphi_i(\bar x) $. Then there exist $x_k\stackrel{\varphi_i}\rightarrow \bar x$ and $x_k^\ast \stackrel{w^\ast }\rightarrow x^\ast $ such that
\begin{equation*}
  x_k^\ast \in\hat\partial \varphi(x_k)=\triangledown f(x_k)^\ast (\partial g_i(f(x_k)))\ \ \forall k.
\end{equation*}
Thus for each $k$, there is $y_k^\ast \in\partial g_i(f(x_k))$ such that 
\begin{equation}\label{4.18}
  x_k^\ast =\triangledown f(x_k)^\ast (y_k^\ast ).
\end{equation}
Note that $g_i$ is a continuous convex function and thus is locally Lipschitz   around $f(\bar x)$. This implies that $\{\|y_k^\ast \|\}$ is bounded. Since $\mathbb{Y}$ is an Asplund space, it follows that $\mathbf{B}_{\mathbb{X}^\ast }$ is weak$^\ast $ sequentially compact and without loss of generality, we can assume that 
$$
y_k^\ast \stackrel{w^\ast }\rightarrow y^\ast \in \mathbb{Y}^\ast 
$$
(taking  a subsequence if necessary). Noting that $y_k^\ast \in\partial g_i(f(x_k))$, $x_k\rightarrow x$ and $\triangledown f$ is continuous at $x$, it follows that
$$
y^\ast \in \partial g_i(f(\bar x))\ \ {\rm and}\ \ \triangledown f(x_k)^\ast (y_k^\ast )\stackrel{w^\ast }\rightarrow \triangledown f(\bar x)^\ast (y^*)
$$
and consequently \eqref{4.18} gives
$$
x^\ast =\triangledown f(\bar x)^\ast (y^\ast )\in\triangledown f(\bar x)^\ast (\partial g_i(f(\bar x)))=\hat\partial\varphi_i(\bar x)
$$
(thanks to $x_k^\ast \stackrel{w^\ast }\rightarrow x^\ast $). Hence \eqref{4.17} holds. The proof is complete.\hfill $\Box$\\[1pt]

\begin{remark}
	
In \eqref{4.12} of \cref{lem4.1}, the first equality means that each $\varphi_i$ is {\it lower regular} at $\bar x$ (see \cite[Definition 1.91]{Mordukhovich}) and the second equality is inspired by \cite[Chain Rule 1]{Clarke1981} that established the calculus of Clarke subdifferentials for convex-composite functions.  

\end{remark}

Now, we are in a position to present dual characterizations for error bounds of the convex-composite   inequality system \eqref{4.7}

\begin{theorem}\label{th4.2}
Consider   the convex-composite inequality system \eqref{4.7}. Then   \eqref{4.7} has a local error bound at $\bar x\in A$ if and only if there exist $\eta,\delta\in (0,+\infty) $   such that for any $x\in A\cap \mathbf{B}(\bar x, \delta)$ and $x^*\in \hat N(A,x)\cap \mathbf{B}_{\mathbb{X}^\ast }$, there are multipliers $\lambda_i\in [0,\eta], i\in I(x)$ satisfying
\begin{equation}\label{4.19}
  x^*\in\sum_{i\in I(x)} \lambda_i \hat\partial\varphi_i(x).
\end{equation}
\end{theorem}

{\bf Proof.} Set  $\varphi_i^+(x):=\max\{\varphi_i(x),0\}$ for any $x$.

{\it The necessity part}. By the local error bound, there exist $\eta,\delta>0$ such that \eqref{4.2} holds. Let $x\in A\cap \mathbf{B}(\bar x, \delta)$ and $x^*\in \hat N(A,x)\cap \mathbf{B}_{\mathbb{X}^\ast }$. Then $x^\ast \in \hat\partial \mathbf{d}(\cdot, A)(x)$ and thus for any $\epsilon>0$, there is $\delta_1\in (0, \delta-\|x-\bar x\|)$ such that
\begin{equation*}\label{4.20}
  \langle x^\ast ,u-x\rangle\leq \mathbf{d}(u,A)-\mathbf{d}(x,A)+\eta\epsilon\|u-x\|,\ \ \forall u\in \mathbf{B}(x,\delta_1)
\end{equation*}
and it follows from \eqref{4.2} that 
\begin{equation*}\label{4.20}
  \langle \frac{x^\ast }{\eta},u-x\rangle\leq \psi(u)+\epsilon\|u-x\|,\ \ \forall u\in \mathbf{B}(x,\delta_1),
\end{equation*}
where $\psi(x):=\sum_{i=1}^{m}\varphi_i^+(x)$. This together  with \cite[Theorem 3.36]{Mordukhovich} and \cite[Theorem 3.46]{Mordukhovich} implies that
$$
\frac{x^\ast }{\eta}\in \hat\partial\psi(x)\subseteq\partial\psi(x)\subseteq\sum_{i=1}^{m}\partial\varphi_i^+(x)=\sum_{i\in I(x)}\partial\varphi_i^+(x)\subseteq\sum_{i\in I(x)}[0,1]\partial\varphi_i(x)
$$
(the equality holds by the continuity of $\varphi_i$). Then there are multipliers $\lambda_i\in [0,\eta], i\in I(x)$ such that
$$
x^\ast\in \sum_{i\in I(x)}\lambda_i \partial\varphi_i(x)
$$
This and \cref{lem4.1} imply that \eqref{4.19} holds.

{\it The sufficiency part}.  Note that each $g_i$ is locally Lipschitz   around $f(\bar x)$ and then there exist $\kappa,r_0>0$ such that
\begin{equation}\label{4.20}
  \|g_i(y_1)-g_i(y_2)\|\leq\kappa\|y_1-y_2\|,\ \ \forall y_1,y_2\in \mathbf{B}(f(\bar x), r_0) \ {\rm and} \ \forall i.
\end{equation}
Denote $L:=\|\triangledown f(\bar x)\|+1$. Let $\epsilon>0$ be such that $m\eta\kappa\epsilon<1$. Since $f$ is continuously differentiable at $\bar x$, there exists $r_1\in (0,\delta)$ with $2Lr_1<r_0$ such that 
\begin{equation}\label{4.21}
 \|\triangledown f(x)\|\leq L\ \ {\rm and} \ \ \|f(x)-f(u)-\triangledown f(u)(x-u)\|<\epsilon\|x-u\|,\ \forall x,u\in \mathbf{B}(\bar x,r_1).
\end{equation}
Take $\delta_1:=\frac{r_1}{2}$ and let $x\in \mathbf{B}(\bar x,\delta_1)\backslash A$. Then $\mathbf{d}(x,A)\leq\|x-\bar x\|<\delta_1$. Choose $\beta\in (0,1)$ such that 
$\beta>\max\{\frac{\mathbf{d}(x,A)}{\delta_1}, m\eta\kappa\epsilon\}$. Applying \cref{lem2.3}, there exist $a\in A$ and $a^*\in\hat N(A,a)$ with $\|a^*\|=1$ such that
\begin{equation}\label{4.22}
  \beta\|x-a\|\leq\min\{\mathbf{d}(x, A), \langle a^*, x-a\rangle\}.
\end{equation}
Then 
$$
\|a-\bar x\|\leq\|a-x\|+\|x-\bar x\|<\frac{\mathbf{d}(x,A)}{\beta}+\delta_1<2\delta_1=r_1
$$
and it follows from \eqref{4.19} and \cref{lem4.1} that there exist $\lambda_i\in[0,\eta],y_i^*\in\partial g_i(f(a)) (\forall i\in I(a))$ such that
\begin{equation*}\label{4.23}
  a^*=\sum_{i\in I(a)}\lambda_i\triangledown f(a)^*(y_i^*).
\end{equation*}
This together with \eqref{4.20} and \eqref{4.21} implies that
\begin{eqnarray*}
\langle a^*, x-a\rangle&=&\sum_{i\in I(a)}\lambda_i\langle\triangledown f(a)^*(y_i^*),x-a\rangle=\sum_{i\in I(a)}\lambda_i\langle y_i^*,\triangledown f(a)(x-a)\rangle\\
&\leq&\sum_{i\in I(a)}\lambda_i\big(g_i(f(a)+\triangledown f(a)(x-a))-g_i(f(a))\big)\\
&\leq&\sum_{i\in I(a)}\lambda_i\big(g_i(f(x))-g_i(f(a))\big)+\sum_{i\in I(a)}\lambda_i\big(g_i(f(a)+\triangledown f(a)(x-a))-g_i(f(x))\big)\\
&\leq&\sum_{i\in I(a)}\lambda_i\varphi_i(x)+\sum_{i\in I(a)}\lambda_i\kappa\|f(x)-f(a)-\triangledown f(a)(x-a))\|\\
&\leq&\eta\sum_{i=1}^{m}\varphi_i^+(x)+m\eta\kappa\epsilon\|x-a\|
\end{eqnarray*}
and it follows from \eqref{4.22} that 
$$
(\beta-m\eta\kappa\epsilon)\|x-a\|\leq\eta\sum_{i=1}^{m}\varphi_i^+(x).
$$
By taking limits as $\beta\uparrow 1$, one gets
$$
\mathbf{d}(x, A)\leq\frac{\eta}{1-m\eta\kappa\epsilon}\sum_{i=1}^{m}\varphi_i^+(x).
$$
This means that the convex-composite  inequality  system \eqref{4.7} has a local error bound with constant $\frac{\eta}{1-m\eta\kappa\epsilon}$. The proof is complete.\hfill$\Box$\\[1pt]

Note that  the convex-composite inequality  system \eqref{4.7} has a global error bound if and only if  the convex-composite inequality system \eqref{4.7} has a local error bound at all points in $A$ with the same constant and so the following corollary on the global error bound of \eqref{4.7} follows from \cref{th4.2} and its proof.

\begin{corollary}\label{coro4.1}
Let the convex-composite   inequality system be given as in \eqref{4.7} and assume that each $g_i$ is Lipschitz. Then  the convex-composite nequality  system \eqref{4.7} has a global error bound if and only if there exists $\eta\in (0,+\infty)$ such that for any $x\in A$ and $x^*\in \hat N(A,x)\cap \mathbf{B}_{\mathbb{X}^\ast }$, there are multipliers $\lambda_i\in [0,\eta], i\in I(x)$ satisfying \eqref{4.19}.
\end{corollary}

For the case of convex  inequality systems defined by convex and continuous functions $g_1,\cdots,g_m$, the following corollary follows from \cref{th4.2} and its proof.

\begin{corollary}\label{coro4.2}
Consider the following convex inequality  system:
  \begin{equation}\label{4.25}
    g_1(y)\leq 0,\cdots,g_m(y)\leq 0.
  \end{equation}
We  denote by $C:=\bigcap_{i=1}^mC_i$ the solution set, where $C_i:=\, \{y\in \mathbb{Y} : g_i(y)\leq 0\}, i=1,\cdots,m$.  Then
\begin{itemize}
  \item [\rm (i)] the  convex  inequality  system \eqref{4.25} has a local error bound at $\bar y\in C$ if and only if there exist $\eta,\delta\in (0,+\infty)$  such that for any $y\in C\cap \mathbf{B}(\bar y, \delta)$ and $y^*\in N(C, y)\cap \mathbf{B}_{\mathbb{Y}^\ast }$, there are multipliers $\alpha_i\in [0,\eta], i\in J(y)$ satisfying
  \begin{equation}\label{4.26}
    y^*\in  \sum_{i\in J(y)}    \alpha_i\partial g_i(y);
  \end{equation}
  \item[\rm (ii)]   the convex inequality system \eqref{4.25} has a global error bound if and only if there exists $\eta\in (0,+\infty)$   such that for any $y\in C$ and $y^*\in N(C, y)\cap \mathbf{B}_{\mathbb{Y}^\ast }$, there are multipliers $\alpha_i\in [0,\eta], i\in J(y)$ satisfying \eqref{4.26},   where $J(y):=\{i:g_i(y)=0\}$.

\end{itemize}
\end{corollary}

\begin{remark}
It is noted that \cite[Theorem 2]{AC1988} and \cite[Corollary 8]{BBL}  provide sufficient conditions for the global error bound of   convex inequality  systems. Compared with these two results,   \cref{coro4.1} actually gives a necessary and sufficient condition for global error bound, which improves results of \cite[Theorem 2]{AC1988} and \cite[Corollary 8]{BBL}.
\end{remark}

\begin{theorem}\label{th4.3}
Let $\bar x\in A$ be such that $\triangledown f(\bar x)$ is surjective. Then convex-composite inequality  system \eqref{4.7} has a local error bound at $\bar x$ if and only if  the convex inequality  system \eqref{4.25} has a local error bound at $f(\bar x)$.
\end{theorem}

{\bf Proof.} Since $\triangledown f(\bar x)$ is surjective, then \cite[Theorem 1.57]{Mordukhovich} implies that  there exist $\kappa,r_0>0$ such that
\begin{equation}\label{4.28}
  \mathbf{d}(x, f^{-1}(y))\leq\kappa \|y-f(x)\|\ \ \forall (x,y)\in \mathbf{B}(\bar x,r_0)\times \mathbf{B}(f(\bar x), r_0).
\end{equation}
Using the proof of \cref{lem3.1}, there are $\ell,L>0$ and $r_1\in (0, r_0)$ such that 
\begin{equation}\label{4.29}
\ell\mathbf{B}_{\mathbb{Y}}\subseteq \triangledown f(x)(\mathbf{B}_{\mathbb{X}}) \ \ {\rm and} \ \  \|\triangledown f(x)\|\leq L,\ \forall x\in \mathbf{B}(\bar x,r_1)
\end{equation}
and
\begin{equation}\label{4.30}
\hat N(A,x)\cap \ell\mathbf{B}_{\mathbb{X}^*}\subseteq \triangledown f(x)^*(N(C,f(x))\cap \mathbf{B}_{\mathbb{Y}^*})\subseteq \hat N(A,x)\cap L\mathbf{B}_{\mathbb{X}^*}, \ \forall x\in \mathbf{B}(\bar x,r_1)\cap A. 
\end{equation}

{\it The necessity part.} By \cref{th4.2}, there exist $\eta,\delta\in (0,+\infty)$  such that for any $x\in A \cap \mathbf{B}(\bar x,\delta)$ and any $x^*\in 	\hat N(A,x)\cap \mathbf{B}_{\mathbb{X}^*}$ there are multipliers $\lambda_i\in [0,\eta], i\in I(x)$ satisfying
\begin{equation*} \label{4.28a}    
x^*\in  \sum_{i\in I(x)} \lambda_i  \hat\partial \varphi_i(x).   
\end{equation*}

Take $\delta_1\in (0,r_1)$ such that $\kappa\delta_1<\delta$. Let $y\in C\cap \mathbf{B}(f(\bar x),\delta_1)$ and $y^*\in N(C,y)\cap \mathbf{B}_{\mathbb{Y}^*}$. Then by \eqref{4.28}, one has
$$
\mathbf{d}(\bar x, f^{-1}(y))\leq\kappa\|f(\bar x)-y\|<\kappa\delta_1
$$
and thus there is $x\in f^{-1}(y)$ such that $\|x-\bar x\|<\kappa\delta_1$. This implies that $x\in A\cap \mathbf{B}(\bar x,\delta)$ and thus $y^*\in N(C,f(x))\cap \mathbf{B}_{\mathbb{Y}^*}$. Note that $I(x)=J(f(x))$ and it follows from \eqref{4.29},  \eqref{4.30} and \cite[Corollary 1.15]{Mordukhovich} that
\begin{eqnarray*}
\triangledown f(x)^*(y^*)\in \triangledown f(x)^*(N(C,f(x)))\cap L\mathbf{B}_{\mathbb{X}^*}=\hat N(A,x)\cap L\mathbf{B}_{\mathbb{X}^*}.
\end{eqnarray*}
Then there are multipliers $\lambda_i\in [0,\eta], i\in I(x)$ such that
\begin{eqnarray*}
\triangledown f(x)^*\Big(\frac{y^*}{L}\Big)\in \sum_{i\in I(x)}  \lambda_i  \hat\partial\varphi_i(x)
&=&  \sum_{i\in I(x)}    \lambda_i   \triangledown f(x)^*(\partial g_i(f(x)))\\
&=&\triangledown f(x)^*\Big(  \sum_{i\in J(f(x))}   \lambda_i   \partial g_i(f(x))\Big)
\end{eqnarray*}
(thanks to $I(x)= J(f(x))$). Since $\triangledown f(x)^*$ is one-to-one, it follows that
$$
y^*\in \sum_{i\in J(y)} L\lambda_i   \partial g_i(f(x)).
$$
This and \cref{coro4.2} imply that  the convex  inequality system \eqref{4.25} has a local error bound at $f(\bar x)$.

{\it The sufficiency part}. By \cref{coro4.2}, there exist $\eta,\delta>0$ such that for any $y\in C\cap \mathbf{B}(f(\bar x), \delta)$ and $y^*\in N(C, y)\cap \mathbf{B}_{\mathbb{Y}^\ast }$, there are multipliers $\alpha_i\in [0,\eta], i\in J(y)$ such that \eqref{4.26} holds. 

Take $\delta_1>0$ such that $f(\mathbf{B}(\bar x,\delta_1))\subseteq \mathbf{B}(f(\bar x),\delta)$. Let $x\in A\cap \mathbf{B}(\bar x,\delta_1)$ and $x^*\in\hat N(A,x)\cap \mathbf{B}_{\mathbb{X}^*}$. By virtue of \eqref{4.30}, there is $y^*\in N(C,f(x))\cap \mathbf{B}_{\mathbb{Y}^*}$ such that
\begin{equation}\label{4.31}
  \ell x^*=\triangledown f(x)^*(y^*).
\end{equation}
From $y^*\in N(C,f(x))\cap \mathbf{B}_{\mathbb{Y}^*}$, there are multipliers $\alpha_i\in [0,\eta], i\in J(y)$ such that \eqref{4.26} holds and it follows from \cref{lem4.1} that
\begin{eqnarray*}
\triangledown f(x)^*(y^*)&\in& \triangledown f(x)^*\Big(\sum_{i\in J(f(x))}\alpha_i\partial g_i(f(x))\Big)\\
 &=& \sum_{i\in J(f(x))} \alpha_i \triangledown f(x)^*\Big(\partial g_i(f(x))\Big)\\
 &=& \sum_{i\in I(x)} \alpha_i\hat\partial \varphi_i(x)
\end{eqnarray*}
(thanks to $I(x)=J(f(x))$). Then \eqref{4.31} gives that 
$$
x^*\in  \sum_{i\in I(x)}\frac{\alpha_i}{\ell}\hat\partial \varphi_i(x).
$$
This and \cref{coro4.1} imply that convex-composite inequality  system \eqref{4.7} has a local error bound at $\bar x$. The proof is complete.\hfill$\Box$

\section{Conclusions}

The study on subtransversality and strong CHIP refers to the so-called convex feasibility problem of finding a common point in the intersection of convex sets in a Hilbert space. Since the two   concepts are proved to ensure a better norm convergence of the sequence produced by projection-based algorithms for solving the convex feasibility problem, subtransversality and strong CHIP, especially relationships with other qualifications in approximation theory and optimization, have been well-recognized and extensively studied. One of  the most important  results on these concepts is that the dual characterization of subtransversality for a collection of closed convex sets can be established via the  strong CHIP and property (G). The work in this paper was  to study such issue for a collection of special closed sets in  Asplund spaces. Although the dual characterization result on subtransversality for general closed sets may not hold (if dropping the convexity assumption), several necessary conditions for subtransversality are still obtained in terms of Fr\'echet/limiting normal cones. Further, the equivalence result among subtransversality, strong Fr\'echet CHIP and property (G) for some special closed sets in convex-composite optimization is proved. Our work is an extension of  a duality characterization of subtransversality via the  strong CHIP and property (G) to the possible non-convex case.

\medskip

\noindent{\bf Acknowlegements} The authors would like to thank the two reviewers for their  careful reading of the paper and for their very insightful comments on the original version of this paper.

\medskip

\noindent{\bf Data Availability} Data sharing not applicable to this article as no datasets were generated or analysed during
the current study.

\medskip

\noindent{\large \bf Declarations}

\medskip

\noindent{\bf Competing Interests}  The authors declare no competing interests.

\bibliography{WTY1}

\begin{thebibliography}{10}

\bibitem{AC1988}
A.~Auslender and J.-P. Crouzeix.
\newblock Global regularity theorems.
\newblock {\em Math. Oper. Res.}, 13(2):243--253, 1988.

\bibitem{AusDanThi05}
D.~Aussel, A.~Daniilidis, and L.~Thibault.
\newblock Subsmooth sets: functional characterizations and related concepts.
\newblock {\em Trans. Amer. Math. Soc.}, 357(4):1275--1301, 2005.

\bibitem{BDL2005}
A.~Bakan, F.~Deutsch, and W.~Li.
\newblock Strong {CHIP}, normality, and linear regularity of convex sets.
\newblock {\em Trans. Amer. Math. Soc.}, 357(10):3831--3863, 2005.

\bibitem{BB}
H.~H. Bauschke and J.~M. Borwein.
\newblock On the convergence of von neumann's alternating projection algorithm
  for two sets.
\newblock {\em Set-Valued Analysis}, 1(2):185--212, 1993.

\bibitem{BB1}
H.~H. Bauschke and J.~M. Borwein.
\newblock On projection algorithms for solving convex feasibility problems.
\newblock {\em SIAM Rev.}, 38(3):367--426, 1996.

\bibitem{BBL97}
H.~H. Bauschke and J.~M. Borwein.
\newblock The method of cyclic projections for convex sets in hilbert space.
\newblock In {\em Recent Developments in Optimization Theory and Nonlinear
  Analysis}, volume 204 of {\em Y. Censor and S. Reich (Eds), Contemporary
  Mathematics}, pages 1--38. Vanderbilt Univ. Press, Nashville, TN, 1997.

\bibitem{BBL}
H.~H. Bauschke, J.~M. Borwein, and W.~Li.
\newblock Strong conical hull intersection property, bounded linear regularity,
  {J}ameson's property {$(G)$}, and error bounds in convex optimization.
\newblock {\em Math. Program., Ser. A}, 86(1):135--160, 1999.

\bibitem{BBT}
H.~H. Bauschke, J.~M. Borwein, and Paul Tseng.
\newblock Bounded linear regularity, strong {CHIP}, and {CHIP} are distinct
  properties.
\newblock {\em J. Convex Anal.}, 7(2):395--412, 2000.

\bibitem{BZ2005}
J.~M. Borwein and Q.J. Zhu.
\newblock {\em Techniques of Variational Analysis}.
\newblock Springer, New York, 2005.

\bibitem{BCK2020}
H.~T. Bui, N.~D. Cuong, and A.~Y. Kruger.
\newblock Transversality of collections of sets: Geometric and metric
  characterizations.
\newblock {\em Vietnam J. Math.}, 48(2):277–297, 2020.

\bibitem{CDW1}
C.~K. Chui, F.~Deutsch, and J.~D. Ward.
\newblock Constrained best approximation in {H}ilbert space.
\newblock {\em Constr. Approx.}, 6(1):35--64, 1990.

\bibitem{CDW2}
C.~K. Chui, F.~Deutsch, and J.D. Ward.
\newblock Constrained best approximation in {H}ilbert space. {II}.
\newblock {\em J. Approx. Theory}, 71(2):213--238, 1992.

\bibitem{Clarke1981}
F.~H. Clarke.
\newblock Generalized gradients of lipschitz functionals.
\newblock {\em Advances in Mathematics}, 40:52--67, 1981.

\bibitem{Combettes96}
P.~L. Combettes.
\newblock The convex feasibility problem in image recovery.
\newblock In P.~Hawkes, editor, {\em Advances in Imaging and Electron Physics},
  volume~95, pages 155--270. Academic Press, New York, 1996.

\bibitem{D}
F.~Deutsch.
\newblock The role of the strong conical hull intersection property in convex
  optimization and approximation.
\newblock In {\em Approximation theory {IX}, {V}ol. {I}. ({N}ashville, {TN},
  1998)}, Innov. Appl. Math., pages 105--112. Vanderbilt Univ. Press,
  Nashville, TN, 1998.

\bibitem{Deutsch98}
F.~Deutsch.
\newblock The role of the strong conical hull intersection property in convex
  optimization and approximation.
\newblock In {\em Approximation theory {IX}, {V}ol. {I}. ({N}ashville, {TN},
  1998)}, Innov. Appl. Math., pages 105--112. Vanderbilt Univ. Press,
  Nashville, TN, 1998.

\bibitem{DLS}
F.~Deutsch, W.~Li, and J.~Swetits.
\newblock Fenchel duality and the strong conical hull intersection property.
\newblock {\em J. Optim. Theory Appl.}, 102(3):681--695, 1999.

\bibitem{DeuLiWar99}
F.~Deutsch, W.~Li, and J.~D. Ward.
\newblock Best approximation from the intersection of a closed convex set and a
  polyhedron in {H}ilbert space, weak {S}later conditions, and the strong
  conical hull intersection property.
\newblock {\em SIAM J. Optim.}, 10(1):252--268, 1999.

\bibitem{DGZ1993}
R.~Deville, G.~Godefroy, and V.~Zizler.
\newblock {\em Smoothness and Renorming in Banach Spaces}.
\newblock Pitman Monogr. Surveys Pure Appl. Math. 64, Longman, Harlow,, 1993.

\bibitem{Dolecki1982}
S.~Dolecki.
\newblock Tangency and differentiation, some applications of convergence
  theory.
\newblock {\em Ann. Math. Pura Appl.}, 130:223–255, 1982.

\bibitem{DonRoc09}
A.~L. Dontchev and R.~T. Rockafellar.
\newblock {\em Implicit functions and solution mappings : a view from
  variational analysis}.
\newblock Implicit functions and solution mappings : a view from variational
  analysis, 2013.

\bibitem{ET2007}
E.~Ernst and M.~Th\'era.
\newblock Boundary half-strips and the strong {CHIP}.
\newblock {\em SIAM Journal on Optimization}, 18(3):834--852, 2007.

\bibitem{F}
M.~Fabian.
\newblock Subdifferentiability and trustworthiness in the light of a new
  variational principle of borwein and preiss.
\newblock {\em Charles University in Prague}, 30(2), 1989.

\bibitem{FM1998}
M.~Fabian and B.~S. Mordukhovich.
\newblock Nonsmooth characterizations of {A}splund spaces and smooth
  variational principles.
\newblock {\em Set-Valued Analysis}, 6:381--406, 1998.

\bibitem{HW}
H.~Hu and Q.~Wang.
\newblock Strong {{CHIP}} for infinite systems of convex sets in normed linear
  spaces.
\newblock {\em Optimization}, 59(2):235--251, 2010.

\bibitem{NT2004}
V.~N. Huynh and M.~Th{\'e}ra.
\newblock Error bounds and implicit multifunction theorem in smooth {B}anach
  spaces and applications to optimization.
\newblock {\em Set-Valued Anal.}, 12(1-2):195--223, 2004.

\bibitem{Ioffe1989}
A.~D. Ioffe.
\newblock Approximate subdifferentials and applications 3: the metric theory.
\newblock {\em Mathematika}, 36(1):1--38, 1989.

\bibitem{Ioffe2017}
A.~D. Ioffe.
\newblock {\em Variational Analysis of Regular Mappings: Theory and
  Applications}.
\newblock Springer Monographs in Mathematics. Springer, 2017.

\bibitem{IO2008}
A.~D. Ioffe and J.V. Outrata.
\newblock On metric and calmness qualification conditions in subdifferential
  calculus.
\newblock {\em Set-valued and variational analysis}, 16(2-3):199--227, 2008.

\bibitem{IP1996}
A.~D. Ioffe and J-P. Penot.
\newblock Subdifferentials of performance functions and calculus of
  coderivatives of set-valued mappings.
\newblock {\em Serdica Math J}, 22(22):359--384, 1996.

\bibitem{Jameson}
G.~Jameson.
\newblock {\em Ordered linear spaces}.
\newblock Springer Berlin Heidelberg, 1970.

\bibitem{J}
G.~Jameson.
\newblock The duality of pairs of wedges.
\newblock {\em Proceedings of the London Mathematical Society}, 24(3):531--547,
  1972.

\bibitem{J1}
V.~Jeyakumar.
\newblock The strong conical hull intersection property for convex programming.
\newblock {\em Mathematical Programming}, 106(1):81--92, 2006.

\bibitem{KL1997}
K.~C. Kiwiel and B.~Lopuch.
\newblock Surrogate projection methods for finding fixed points of firmly
  nonexpansive mappings.
\newblock {\em SIAM Journal on Optimization}, 7(4), 1997.

\bibitem{Kru15}
A.~Y. Kruger.
\newblock Error bounds and metric subregularity.
\newblock {\em Optimization: A Journal of Mathematical Programming and
  Operations Research}, 2015.

\bibitem{KLT2017}
A.~Y. Kruger, D.~R Luke, and N.~H. Thao.
\newblock About subtransversality of collections of sets.
\newblock {\em Set-valued and variational analysis}, 25(4):701–729, 2017.

\bibitem{KLT2018}
A.~Y. Kruger, D.~R. Luke, and N.~H. Thao.
\newblock Set regularities and feasibility problems.
\newblock {\em Mathematical Programming}, 168(1-2):279–311, 2018.

\bibitem{LP}
A.~S. Lewis and J.~S. Pang.
\newblock {\em Error Bounds for Convex Inequality Systems}.
\newblock In: Crouzeix, J. P., ed., Generalized Convexity, Proceedings of the
  Fifth Sysposium on Generalized Convexity, Luminy Marseille, 1997, pp. 75-100,
  1997.

\bibitem{LJ}
C.~Li and X.~Q. Jin.
\newblock Nonlinearly constrained best approximation in hilbert spaces: The
  strong {CHIP} and the basic constraint qualification.
\newblock {\em SIAM Journal on Optimization}, 13(1):228--239, 2002.

\bibitem{LN}
C.~Li and K.~F. Ng.
\newblock Constraint qualification, the strong {CHIP}, and best approximation
  with convex constraints in banach spaces.
\newblock {\em SIAM Journal on Optimization}, 14(2):584--607, 2003.

\bibitem{LN1}
C.~Li and K.~F. Ng.
\newblock Strong {CHIP} for infinite system of closed convex sets in normed
  linear spaces.
\newblock {\em SIAM Journal on Optimization}, 16(2):311--340, 2005.

\bibitem{LN2}
C.~Li and K.~F. Ng.
\newblock The dual normal {CHIP} and linear regularity for infinite systems of
  convex sets in banach spaces.
\newblock {\em SIAM Journal on Optimization}, 24(3):1075--1101, 2014.

\bibitem{LNP}
C.~Li, K.~F. Ng, and T~K Pong.
\newblock The secq, linear regularity, and the strong {CHIP} for an infinite
  system of closed convex sets in normed linear spaces.
\newblock {\em SIAM Journal on Optimization}, 18:643--665, 2007.

\bibitem{Mordukhovich}
B.~S. Mordukhovich.
\newblock {\em Variational Analysis and Generalized Differentiation I}.
\newblock Springer, New York, NY, 2006.

\bibitem{MS}
B.~S. Mordukhovich and Y.~Shao.
\newblock Nonsmooth sequential analysis in {A}splund spaces.
\newblock {\em Transactions of the American Mathematical Society},
  348:1235--1280, 1996.

\bibitem{MW2000}
B.~S. Mordukhovich and B.~Wang.
\newblock On variational characterizations of {A}splund spaces.
\newblock In {\em Constructive, experimental, and nonlinear analysis
  ({L}imoges, 1999)}, volume~27 of {\em CRC Math. Model. Ser.}, pages 245--254.
  CRC, Boca Raton, FL, 2000.

\bibitem{moreau}
J.~J. Moreau.
\newblock Fonctionnelles convexes.
\newblock {\em S\'eminaire Jean Leray}, (2):1--108, 1966-1967.

\bibitem{NY}
K.~F Ng and W.~H. Yang.
\newblock Error bounds for abstract linear inequality systems.
\newblock {\em SIAM Journal on Optimization}, 13(1):24--43, 2002.

\bibitem{NT2001}
H.~V. Ngai and M.~Th\'era.
\newblock Metric inequality, subdifferential calculus and applications.
\newblock {\em Set Valued Analysis}, 9(1-2):187--216, 2001.

\bibitem{P}
J.~S. Pang.
\newblock Error bounds in mathematical programming.
\newblock {\em Mathematical Programming}, 79(1):299--332, 1997.

\bibitem{Penot2013}
J.~P. Penot.
\newblock {\em Calculus without Derivatives}.
\newblock 2013.

\bibitem{Phelps}
R.~R. Phelps.
\newblock {\em Convex Functions, Monotone Operators and Differentiability}.
\newblock Lecture Notes in Math. 1364, Springer, New York, 1989.

\bibitem{Ra}
H.~R{\aa}dstr\"{o}m.
\newblock An embedding theorem for spaces of convex sets.
\newblock {\em Proceedings of the American Mathematical Society},
  3(1):165--169, 1952.

\bibitem{Rudin}
W.~Rudin.
\newblock {\em Functional analysisd}.
\newblock McGraw-Hill, New York, 1991.

\bibitem{thibault}
L.~Thibault.
\newblock {\em Unilateral Variational Analysis in Banach Spaces}.
\newblock World Scientific, 2023.

\bibitem{WH}
Z.~Wei and Q.~He.
\newblock Normal property, jameson property, {CHIP} and linear regularity for
  an infinite system of convex sets in banach spaces.
\newblock {\em Optimization: A Journal of Mathematical Programming and
  Operations Research}, 65(12):2095--2114, 2016.

\bibitem{WZ}
Z.~Wei, J.-C. Yao, and X.~Y. Zheng.
\newblock Strong {A}badie {CQ}, {ACQ}, calmness and linear regularity.
\newblock {\em Mathematical Programming}, 145:97--131, 2014.

\bibitem{W}
Zhou Wei.
\newblock Linear regularity for an infinite system formed by {$p$}-uniformly
  subsmooth sets in {B}anach spaces.
\newblock {\em Taiwanese J. Math.}, 16(1):335--352, 2012.

\bibitem{Yost93}
D.~Yost.
\newblock Asplund spaces for beginners.
\newblock {\em Acta Univ. Carolin. Math. Phys.}, 34(2):159--177, 1993.
\newblock Selected papers from the 21st Winter School on Abstract Analysis
  (Pod{\v{e}}brady, 1993).

\bibitem{ZN2004}
X.~Y. Zheng and K.~F. Ng.
\newblock Metric regularity and constraint qualifications for convex
  inequalities on banach spaces.
\newblock {\em SIAM Journal on Optimization}, 14(3):757--772, 2004.

\bibitem{ZN}
X.~Y. Zheng and K.~F. Ng.
\newblock Linear regularity for a collection of subsmooth sets in banach
  spaces.
\newblock {\em SIAM Journal on Optimization}, 19(1):62--76, 2008.

\bibitem{ZW}
X.~Y. Zheng, Z.~Wei, and J.-C. Yao.
\newblock Uniform subsmoothness and linear regularity for a collection of
  infinitely many closed sets.
\newblock {\em Nonlinear Analysis: Theory}, 73(2):413--430, 2010.

\end{thebibliography}

\bibliographystyle{plain}

\end{document}